   \documentstyle[11pt,leqno]{article}
   \newtheorem{dfn}{Definition}[section]
   \newtheorem{lem}[dfn]{Lemma}
   \newtheorem{pro}[dfn]{Proposition}
   \newtheorem{thm}[dfn]{Theorem}
   \newtheorem{rem}[dfn]{Remark}
   \newtheorem{cor}[dfn]{Corollary}
   \newtheorem{exa}[dfn]{Example}
   \newtheorem{cla}{Claim}
   \newtheorem{ste}{Step}
   \def\longhookrightarrow{\lhook\joinrel\longrightarrow}
   
   \def\Biguparrow{\vphantom{\bigg|}\Big\uparrow}
   \def\Bigdownarrow{\vphantom{\bigg|}\Big\downarrow}

   \def\cart{\lower2pt\hbox{\vbox{\hrule width9pt height0.4pt
             \hbox to9pt{\vrule height9pt width 0.4pt
             \hfil\vrule height9pt width 0.4pt}
             \hrule width9pt height 0.4pt}}}
   \def\Spec{\mathop{\rm Spec}\nolimits}
   \def\Ker{\mathop{\rm Ker}\nolimits}
   \def\Coker{\mathop{\rm Coker}\nolimits}
   \def\Hom{\mathop{\rm Hom}\nolimits}
   \def\Homs{\mathop{{\cal H}om}\nolimits}
   \def\Extr{\mathop{\rm Ext}\nolimits}
   
   \def\Ders{\mathop{{\cal D}er}\nolimits}
   \def\Obj{\mathop{\rm Obj}\nolimits}
   \def\gp#1{#1^{\rm gp}}
   \def\sat#1{#1^{\rm sat}}
   \def\tor#1{#1_{\rm tor}}
   \def\fr#1{#1_{\rm f}}
   \def\Dlog{D{\rm log}}
   \def\dlog{d{\rm log}}
   \newcommand{\Ens}{\mbox{\bf Ens}}
   \newcommand{\LS}{{\bf LSch}}
   \newcommand{\LSf}{{\bf LSch}^{\rm f}}
   \newcommand{\LSfs}{{\bf LSch}^{\rm fs}}
   \def\underrel#1#2{\mathrel{\mathop{#1}\limits_{#2}}}
   \newcommand{\pf}{\noindent{\sc Proof.}\hspace{2mm}}
   \newcommand{\qed}{$\Box$}
   \newcommand{\Na}{{\bf N}}
   \newcommand{\Z}{{\bf Z}}
   
   \newcommand{\R}{{\bf R}}
   
   \newcommand{\D}{{\bf LD}}
   \newcommand{\I}{{\cal I}}
   \renewcommand{\L}{{\cal L}}
   \newcommand{\M}{{\cal M}}
   \newcommand{\Mt}{\widetilde{\M}}
   \newcommand{\N}{{\cal N}}
   \renewcommand{\P}{{\cal P}}
   \renewcommand{\O}{{\cal O}}
   \newcommand{\Oi}{{\cal O}^\times}
   \newcommand{\T}{{\cal T}}
   \newcommand{\U}{{\cal U}}
   \newcommand{\X}{\underline{X}}
   \newcommand{\Xp}{\underline{X}'}
   \newcommand{\Xt}{\widetilde{X}}
   \newcommand{\Xtd}{\underline{\widetilde{X}}}
   \newcommand{\Y}{\underline{Y}}
   \newcommand{\Zu}{\underline{Z}}
   \newcommand{\k}{\underline{k}}
   \newcommand{\A}{\underline{A}}
   \newcommand{\uu}{\underline{u}}
   \newcommand{\vv}{\underline{v}}
   \newcommand{\ft}{\widetilde{f}}
   \newcommand{\la}{\lambda}
   \newcommand{\La}{\Lambda}
   \newcommand{\Ext}{{\cal E}xt}
   \newcommand{\Img}{\mbox{\rm Image}\,}
   \newcommand{\art}{{\cal C}_{\Lambda[[Q]]}}
   \newcommand{\Art}{\widehat{{\cal C}}_{\Lambda[[Q]]}}
   \newcommand{\ssd}[1]{\scriptscriptstyle{(#1)}}
   \newcommand{\cvr}{\U=\{X_\la\hookrightarrow V_\la\}}
   \newcommand{\zl}{z^{\ssd{\la}}}
   \newcommand{\zm}{z^{\ssd{\mu}}}
   \newcommand{\zn}{z^{\ssd{\nu}}}
   \newcommand{\ztl}{\zeta^{\ssd{\la}}}
   \newcommand{\ztm}{\zeta^{\ssd{\mu}}}
   \newcommand{\ztn}{\zeta^{\ssd{\nu}}}
   \newcommand{\ulm}{u^{\ssd{\la\mu}}}
   
   \newcommand{\umn}{u^{\ssd{\mu\nu}}}

   \newcommand{\uln}{u^{\ssd{\la\nu}}}
   \newcommand{\slm}{\sigma^{\ssd{\la\mu}}}
   
   \newcommand{\smn}{\sigma^{\ssd{\mu\nu}}}

   \newcommand{\sln}{\sigma^{\ssd{\la\nu}}}
   \title{\bf Log Smooth Deformation Theory}
   \author{{\sc Fumiharu Kato} \\
    Department of Mathematics \\
    Faculty of Science \\
    Kyoto University \\
    Kyoto 606--01 Japan \\
    e--mail: kato@kusm.kyoto--u.ac.jp \thanks{
   \hspace*{1.5em}{\em $1991$ Mathematics Subject Classification\/}.
   Primary 13D10;
   Secondary 14D15, 14M25, 16S80.
   }}
   \date{}
   
\begin{document}
   \def\footnotemark{\relax}
   \maketitle
   \begin{abstract}
   This paper gives a foundation of log smooth deformation theory.
   We study the infinitesimal liftings of log smooth morphisms and show
   that the log smooth deformation functor has a representable hull.
   This deformation theory gives, for example, the following
   two types of deformations: (1) relative deformations of a certain
   kind of a pair of an algebraic variety and a divisor of it, and
   (2) global smoothings of normal crossing varieties.
   The former is a generalization of the relative deformation theory
   introduced by Makio,
   and the latter coincides with the logarithmic deformation theory
   introduced by Kawamata and Namikawa.
   \end{abstract}
   \section{Introduction}
  In this article, we formulate and develop the theory of
{\it log smooth deformation}. Here, log smoothness (more precisely,
logarithmic smoothness) is a concept in {\it log geometry} which is
a generalization of ``usual'' smoothness of morphisms of algebraic
varieties. Log geometry is a beautiful geometric theory which
succesfully generalizes and unifies the scheme theory and the theory
of torus embeddings. This theory was first planned and founded by
Fontaine and Illusie, based on their idea of {\it log structures} on
schemes, and further developed by Kato \cite {Kat1}.
Recently, the importance of log
geometry comes to be recognized by many geometers and applied to
various fields of algebraic and arithmetic geometry. One of such
applications can be seen in the recent work of Steenbrink \cite {Ste1}.
In the present paper, we attempt to apply log geometry to extend the usual
smooth deformation theory by using the concept of log smoothness.

  Log smoothness is one of the most important concepts in log geometry,
and is a log geometric generalization of usual smoothness. For example,
varieties with toric singularities or normal crossing varieties may
become log smooth over certain logarithmic points. Kato \cite {Kat1} showed
that any log smooth morphism is written \'{e}tale locally by the composition
of a usual smooth morphism and a morphism induced by a homomorphism
of monoids which essentially determines the log structures
(Theorem \ref{lisse}). On the other hand, log
smoothness is described by means of {\it log differentials} and
{\it log derivations} similarly to usual smoothness by means of
differentials and derivations. Hence if we consider the log smooth
deformation by analogy with the usual smooth deformation, it is
expected that
the first order deformation is controled by the sheaf of
log derivations. This intuition motivated this work and we shall see
later that this is, in fact, the case.

  In the present paper, we construct
log smooth deformation functor by the concept of infinitesimal log smooth
lifting. The goal of this paper is to show that this functor has a
representable hull in the sense of Schlessinger
\cite {Sch1}, under certain conditions on the underlying schemes
(Theorem \ref{hull}). At the end of this paper,
we give two examples of our log smooth
deformation theory, which are summarized as follows:

\vspace{3mm}
1. {\sc Deformations with divisors} (\S \ref{exam1}): Let $X$ be a
variety over a field $k$. Assume that the variety $X$ is covered
by \'{e}tale open sets
which are smooth over affine torus embeddings, and there exists a divisor
$D$ of $X$ which is the union of the closures of
codimension 1 torus orbits. Then, there exists a log
structure $\M$ on $X$ such that the log scheme $(X, \M)$ is log smooth
over $k$ with trivial
log structure. (The converse is also true in a certain excellent category
of log schemes.)
In this case, our log smooth deformation is a deformation of the
pair $(X,D)$. If $X$ itself is smooth and $D$ is a smooth divisor of $X$,
our deformation coincides with the relative deformation studied by Makio
\cite {Mak1}.

\vspace{3mm}
2. {\sc Smoothings of normal crossing varieties} (\S \ref{exam2}):
If a connected scheme of finite type $X$ over a field $k$ is,
\'{e}tale locally, isomorphic to an affine normal crossing variety
$\Spec k[z_1,\ldots,z_n]/(z_1\cdots z_d)$, then we call $X$ a
normal crossing variety over $k$. If $X$ is $d$--semistable
(cf.\ \cite {Fri1}), there
exists a log structure $\M$ on $X$ such that $(X,\M)$ is log smooth over a
standard log point $(\Spec k,\Na)$ (Theorem \ref{dss}). Then, our log
smooth deformation is nothing but a smoothing of $X$. If the singular
locus of $X$ is connected, our deformation theory coincides
with the one introduced by Kawamata and Namikawa \cite {K-N1}.

\vspace{3mm}
The composition of this paper is as follows. We recall some basic notions of
log geometry in the next section, and review the definition and basic
properties of log smoothness in section 3. In section 4, we study the
characterization of log smoothness by means of the theory of torus
embeddings according to Illusie \cite {Ill1} and Kato \cite {Kat1}. In
section 5, we recall the definitions and basic properties of log derivations
and log differentials. In section 6 and section 7, we give proofs of
theorems stated in section 4. Section 8 is devoted to the formulation of
log smooth deformation theory. This section is the main section of this
present paper. We prove the existance of a representable hull of the
log smooth deformation functor in section 9. In section 10 and section 11,
we give two examples of log smooth deformation. For the reader's convenience,
in section 12, we give a proof of the result of Kawamata and Namikawa
\cite {K-N1} which is relevant to our log smooth deformation.

The author would like to express his thanks to Professors Kazuya Kato and
Yoshinori Namikawa for
valuable suggestions and advice. The author is also very grateful to
Professor Luc Illusie for valuable advice on this paper.

\vspace{3mm}
{\sc Convention}. We assume that all monoids are commutative and have
neutral elements. A homomorphism of monoids is assumed to preserve
neutral elements.
We write the binary operations of all monoids
multiplicatively except in the case of $\Na$ (the monoid of non--negative
integers), $\Z$, etc., which
we write additively.
All sheaves on schemes are considered
with respect to the \'{e}tale topology.
   \section{Fine saturated log schemes}
In this and subsequent sections, we use the terminology of log geometry
basically as in \cite {Kat1}.
Let $X$ be a scheme.
We view the structure sheaf $\O_X$ of $X$ as a sheaf of monoids under
multiplication.

\begin{dfn}{\rm (cf.\ \cite [\S 1]{Kat1})
A {\it pre--log structure} on $X$ is a homomorphism $\M\rightarrow\O_X$ of
sheaves of monoids where $\M$ is a sheaf of monoids on $X$.
A pre--log structure $\alpha:\M\rightarrow\O_X$ is said to be a
{\it log structure} on $X$ if $\alpha$ induces an isomorphism
   $$
   \alpha:\alpha^{-1}(\Oi_X)\stackrel{\sim}{\longrightarrow}\Oi_X.
   $$}
\end{dfn}

\noindent
Given a pre--log structure $\alpha:\M\rightarrow\O_X$, we can
construct the {\it associated log structure}
$\alpha^{\rm a}:\M^{\rm a}\rightarrow\O_X$
functorially by
   \begin{equation}\label{asslog1}
   \M^{\rm a}=(\M\oplus\Oi_X)/\P
   \end{equation}
and
   $$
   \alpha^{\rm a}(x,u)=u\cdot\alpha(x)
   $$
for $(x,u)\in\M^{\rm a}$, where $\P$ is the submonoid defined by
   $$
   \P=\{(x,\alpha(x)^{-1})\: |\: x\in\alpha^{-1}(\Oi_X)\}.
   $$
Here, in general, the quotient $M/P$ of a monoid $M$ with respect
to a submonoid
$P$ is the coset space $M/\sim$ with induced monoid structure,
where the equivalence relation $\sim$ is defined by
   $$
   x\sim y\Leftrightarrow xp=yq \,\mbox{ for some $p,q\in P$}.
   $$
$\M^{\rm a}$ has a universal mapping property:
if $\beta:\N\rightarrow\O_X$ is a log structure on $X$ and
$\varphi:\M\rightarrow\N$ is a homomorphism of sheaves of monoids such that
$\alpha=\beta\circ\varphi$, then there exists a unique lifting
$\varphi^{\rm a}:\M^{\rm a}\rightarrow\N$.
Note that the monoid $\M^{\rm a}$ defined by (\ref{asslog1})
is the push--out of the diagram
   $$
   \M\supset
   \alpha^{-1}(\Oi_X)\stackrel{\alpha}{\longrightarrow}\Oi_X
   $$
in the category of monoids, and the homomorphism $\alpha^{\rm a}$ is
induced by
$\alpha$ and the inclusion $\Oi_X\hookrightarrow\O_X$. We sometimes denote
the monoid $\M^{\rm a}$ by $\M\oplus_{\alpha^{-1}(\Oi_X)}\Oi_X$.
Note that we have the natural isomorphism
   \begin{equation}\label{basic1}
   \M/\alpha^{-1}(\Oi_X)\stackrel{\sim}{\longrightarrow}
   \M^{\rm a}/\Oi_X.
   \end{equation}

\begin{dfn}{\rm
By a {\it log scheme}, we mean a pair $(X,\M)$ with a scheme $X$ and
a log structure $\M$ on $X$. A {\it morphism} of log schemes
$f:(X,\M)\rightarrow(Y,\N)$ is a pair $f=(f,\varphi)$ where
$f:X\rightarrow Y$ is a morphism of schemes and
$\varphi:f^{-1}\N\rightarrow\M$ is a homomorphism of sheaves of monoids
such that the diagram
   $$
   \begin{array}{ccc}
   f^{-1}\N&\stackrel{\varphi}{\longrightarrow}&\M\\
   \Bigdownarrow&&\Bigdownarrow\\
   f^{-1}\O_Y&\longrightarrow&\O_X
   \end{array}
   $$
commutes. }
\end{dfn}

\begin{dfn}\label{logequiv}{\rm
Let $\alpha:\M\rightarrow\O_X$ and $\alpha':\M'\rightarrow\O_X$ be log
structures on a scheme $X$. These log structures are said to be
{\it equivalent} if there exists an isomorphism
$\varphi:\M\stackrel{\sim}{\rightarrow}\M'$ such that
$\alpha=\alpha'\circ\varphi$, i.e., there exists an isomorphism of log schemes
$(X,\M)\stackrel{\sim}{\rightarrow}(X,\M')$ whose underlying morphism of
schemes is the identity ${\rm id}_X$.
Let $\beta:\N\rightarrow\O_Y$ and $\beta':\N'\rightarrow\O_Y$ be log
structures on a scheme $Y$. Let $f:(X,\M)\rightarrow(Y,\N)$ and
$f':(X,\M')\rightarrow(Y,\N')$ be morphisms of log schemes. Then $f$ and $f'$
are said to be {\it equivalent} if there exist isomorphisms
$\varphi:\M\stackrel{\sim}{\rightarrow}\M'$ and
$\psi:\N\stackrel{\sim}{\rightarrow}\N'$ such that
$\alpha=\alpha'\circ\varphi$, $\beta=\beta'\circ\psi$ and the diagram
   $$
   \begin{array}{ccc}
   \M&\stackrel{\varphi}{\longrightarrow}&\M'\\
   \Biguparrow&&\Biguparrow\\
   f^{-1}\N&\underrel{\longrightarrow}{f^{-1}\psi}&f^{-1}\N'
   \end{array}
   $$
commutes. }
\end{dfn}

\noindent
We denote the category of log schemes by $\LS$. For $(S,\L)\in\Obj(\LS)$, we
denote the category of log schemes over $(S,\L)$ by $\LS_{(S,\L)}$.
The following examples play important roles in the sequel.

\begin{exa}\label{trilog}{\rm
On any scheme $X$, we can define a log structure by the inclusion
$\Oi_X\hookrightarrow\O_X$, called the {\it trivial} log structure.
Thus, we have an inclusion functor from the category of schemes to that
of log schemes. We often denote the log scheme $(X, \Oi_X\hookrightarrow\O_X)$
simply by $X$.
}
\end{exa}

\begin{exa}\label{canlog}{\rm
Let $A$ be a commutative ring. For a monoid $P$ , we can define a log structure
canonically on the scheme $\Spec A[P]$, where $A[P]$ denotes the monoid ring
of $P$ over $A$, as the log structure associated to the natural homomorphism,
$$
P\stackrel{\alpha}{\longrightarrow}A[P].
$$
This log structure is called the {\it canonical log structure} on
$\Spec A[P]$.
Thus we obtain a log scheme which we denote simply by $(\Spec A[P], P)$.
A monoid homomorphism $P\rightarrow Q$
induces a morphism $(\Spec A[Q], Q)\rightarrow(\Spec A[P], P)$ of log
schemes. Thus, we have a contravariant functor from the category of monoids
to $\LS_{\Spec A}$.
}
\end {exa}

\begin{exa}\label{torlog}{\rm
Let $\Sigma$ be a fan on $N_\R=\R^d$, $N=\Z^d$, and $X_{\Sigma}$
the toric variety determined by the fan $\Sigma$ over a commutative
ring $A$. Then, we get an induced log structure on the scheme
$X_{\Sigma}$ by gluing the log structures associated to the
homomorphisms
$$
M\cap\sigma^{\vee}\longrightarrow A[M\cap\sigma^{\vee}],
$$
for each cone $\sigma$ in $\Sigma$, where $M=\Hom_{\Z}(N, \Z)$. Thus, a toric
variety $X_{\Sigma}$ is naturally viewed as a log scheme over $\Spec A$,
which we denote by $(X_{\Sigma}, \Sigma)$.
}
\end{exa}

Next, we define important subcategories of $\LS$. These subcategories are
closely related with {\it charts} defined as follows.

\begin{dfn}{\rm
Let $(X,\M)\in\Obj(\LS)$. A {\it chart} of $\M$ is a homomorphism
$P\rightarrow\M$ from the constant sheaf of a monoid $P$ which induces
an isomorphism from the associated log structure $P^{\rm a}$ to $\M$. }
\end{dfn}

\begin{dfn}{\rm
Let $f:(X,\M)\rightarrow(Y,\N)$ be a morphism in $\LS$. A {\it chart} of
$f$ is a triple $(P\rightarrow\M,Q\rightarrow\N,Q\rightarrow P)$, where
$P\rightarrow\M$ and $Q\rightarrow\N$ are charts of $\M$ and $\N$,
respectively, and $Q\rightarrow P$ is a homomorphism for which the diagram
   $$
   \begin{array}{ccc}
   Q&\longrightarrow&P\\
   \Bigdownarrow&&\Bigdownarrow\\
   f^{-1}\N&\longrightarrow&\M
   \end{array}
   $$
is commutative. }
\end{dfn}

\begin{dfn}{\rm (cf.\ \cite [\S 2]{Kat1})
A log structure $\M\rightarrow\O_X$ on a scheme $X$
is said to be {\it fine} if $\M$ has \'{e}tale
locally a chart $P\rightarrow\M$ with $P$ a finitely generated integral
monoid. Here, in general, a monoid $M$ is said to be
{\it finitely generated} if there exists a surjective homomorphism
$\Na^n\rightarrow M$ for some $n$, and a monoid $M$ is
said to be {\it integral} if $M\rightarrow\gp{M}$ is injective, where
$\gp{M}$ denotes the Grothendieck group associated with $M$.
A log scheme $(X,\M)$ with a fine log structure $\M\rightarrow\O_X$ is called
a {\it fine} log scheme. }
\end{dfn}

\noindent
We denote the category of fine log schemes by $\LSf$. Similarly, we denote the
category of fine log schemes over $(S,\L)\in\Obj(\LSf)$ by $\LSf_{(S,\L)}$,
The category $\LSf$ (resp. $\LSf_{(S,\L)}$) is a full subcategory of $\LS$
(resp. $\LS_{(S,\L)}$). Both $\LS$ and $\LSf$ have fiber products.
But the inclusion functor $\LSf\hookrightarrow\LS$ does not preserve fiber
products (cf. Lemma \ref{fpro}). The inclusion functor $\LSf\hookrightarrow
\LS$ has a right adjoint $\LS\rightarrow\LSf$ \cite[(2.7)]{Kat1}. Then, the
fiber product of a diagram $(X,\M)\rightarrow(Z,\P)\leftarrow(Y,\N)$ in
$\LSf$ is the image of that in $\LS$ by this adjoint functor. Note that the
underlying scheme of the fiber product of
$(X,\M)\rightarrow(Z,\P)\leftarrow(Y,\N)$ in $\LS$ is $X\times_{Z}Y$, but
this is not always the case in $\LSf$.

\vspace{3mm}
Next, we introduce more excellent subcategory of $\LS$.

\begin{dfn}\label{satin}{\rm
Let $M$ be a monoid and $P$ a submonoid of $M$.
The monoid $P$ is said to be
{\it saturated} in $M$ if $x\in M$ and $x^n\in P$ for some positive
integer $n$ imply
$x\in P$. An integral monoid $N$ is said to be
{\it saturated} if $N$ is saturated in $\gp{N}$. }
\end{dfn}

\begin{exa}{\rm
Put $M=\Na$ and $P=l\cdot M$ for an integer $l>1$.
Then $P$ is saturated but not
saturated in $M$. }
\end{exa}

\begin{dfn}{\rm
A fine log scheme $(X,\M)\in\Obj (\LSf)$ is said to be {\it saturated} if
the log structure $\M$ is a sheaf of saturated monoids. }
\end{dfn}

\noindent
We denote the category of fine saturated log
schemes by $\LSfs$. Similarly, we denote the category of fine saturated
log schemes over $(S,\L)\in\Obj(\LSfs)$ by $\LSfs_{(S,\L)}$.
The category $\LSfs$ (resp. $\LSfs_{(S,\L)}$) is a full subcategory of $\LSf$
(resp. $\LSf_{(S,\L)}$).
The following lemma is an easy consequence of
\cite [Lemma (2.10)]{Kat1}.

\begin{lem}
Let $f:(X,\M)\rightarrow(Y,\N)$ be a morphism in $\LSfs$, and $Q\rightarrow\N$
a chart of $\N$, where $Q$ is a finitely generated integral saturated
monoid. Then there exists \'{e}tale locally a chart
$(P\rightarrow\M,Q\rightarrow\N,Q\rightarrow P)$
of $f$ extending $Q\rightarrow\N$ such that the monoid $P$ is also finitely
generated, integral and saturated.
\end{lem}

\begin{lem}
The inclusion functor $\LSfs\longhookrightarrow\LSf$ has a right adjoint.
\end{lem}

\pf
Let $M$ be an integral monoid. Define
   $$
   \sat{M}=\{x\in\gp{M}\: |
   \: x^n\in M \;\mbox{for some positive integer $n$}\}.
   $$
 Then $\sat{M}$ is an integral saturated monoid. For any
integral saturated monoid $N$ and homomorphism $M\rightarrow N$, there
exists a unique lifting $\sat{M} \rightarrow N$. In this sense, $\sat{M}$
is the universal saturated monoid associated with $M$. Let $(X,\M)$ be a
fine log scheme. Then we have \'{e}tale locally a chart $P\rightarrow\M$.
This chart defines a morphism $X\rightarrow\Spec \Z[P]$ \'{e}tale locally.
Let $X'=X\times_{\Spec \Z[P]}\Spec \Z[\sat{P}]$. Then $X'\rightarrow
\Spec \Z[\sat{P}]$ induces a log structure $\M'$
by the associated log structure of
$\sat{P}\rightarrow\Z[\sat{P}]\rightarrow\O_{X'}$, and $(X',\M')$ is a fine
saturated log scheme. This procedure defines
a functor $\LSf\rightarrow\LSfs$. It is easy to see that this functor is the
right adjoint of the inclusion functor $\LSfs \hookrightarrow\LSf$. \qed

\begin{cor}
$\LSfs$ has fiber products. More precisely, the fiber product of
morphisms $(X,\M)\rightarrow(Z,\P)\leftarrow(Y,\N)$ in $\LSfs$ is the
image of that in $\LSf$ by the right adjoint functor of
$\LSfs\longhookrightarrow\LSf$.
\end{cor}
   \section{Log smooth morphisms}
In this section, we review the definition and basic properties of log
smoothness (cf. \cite {Kat1}).

\begin{dfn}\label{pback}{\rm
Let $f:X\rightarrow Y$ be a morphism of schemes, and $\N$ a log structure
on $Y$. Then the {\it pull--back} of $\N$, denoted by $f^{*}\N$, is the log
structure on $X$ associated with the pre--log structure
$f^{-1}\N\rightarrow f^{-1}\O_Y\rightarrow\O_X$. A morphism of log schemes
$f:(X,\M)\rightarrow(Y,\N)$ is said to be {\it strict} if the induced
homomorphism $f^{*}\N\rightarrow\M$ is an isomorphism.
A morphism of log schemes
$f:(X,\M)\rightarrow(Y,\N)$ is said to be an {\it exact closed immersion} if
it is strict and
$f:X\rightarrow Y$ is a closed immersion in the ususal sense. }
\end{dfn}

\noindent
Exact closed immersions are stable under base change in $\LSf$
\cite [(4.6)]{Kat1}.

\begin{lem}\label{basic2}
Let $\alpha:\M\rightarrow\O_X$ and $\alpha':\M'\rightarrow\O_X$ be fine log
strctures on a scheme $X$ with a homomorphism $\varphi:\M\rightarrow\M'$ of
monoids such that $\alpha=\alpha'\circ\varphi$. Then, $\varphi$ is an
isomorphism if and only if $\varphi\ {\rm mod}\ \Oi_X:\M/\Oi_X
\rightarrow\M'/\Oi_X$ is an isomorphism.
\end{lem}

\noindent
The proof is straightforward.

\begin{lem}\label{basic3}
Let $f:(X,\M)\rightarrow(Y,\N)$ be a morphism of log schemes. Then, we have
the natural isomorphism
   $$
   f^{-1}(\N/\Oi_Y)\stackrel{\sim}{\longrightarrow}f^{*}\N/\Oi_X.
   $$
In particular, $f$ is strict if and only if the induced morphism
   $$
   f^{-1}(\N/\Oi_Y)\stackrel{\sim}{\longrightarrow}\M/\Oi_X.
   $$
\end{lem}

\pf
The first part is easy to see. For the second part, apply (\ref {basic1}) and
Lemma \ref{basic2}.
\qed

\begin{lem}{\rm (cf. \cite [(1.7)]{Kaj1})}\label{fpro}
Let
   \begin{equation}\label{fpro1}
   (X,\M)\rightarrow(Z,\P)\leftarrow(Y,\N)
   \end{equation}
be morphisms in $\LSfs$.
If $(Y,\N)\rightarrow(Z,\P)$ is strict,
then the fiber product of {\rm (\ref{fpro1})} in $\LSfs$ is isomorphic to
that in $\LS$. In particular, the underlying scheme of the fiber product of
{\rm (\ref{fpro1})} in $\LSfs$ is isomorphic to $X\times_{Z}Y$.
\end{lem}

\pf
Let $P\rightarrow\P$ be a chart of $\P$, where $P$ is a finitely generated
integral saturated monoid. Since $(Y,\N)\rightarrow(Z,\P)$ is strict,
$\rightarrow\P\rightarrow\N$ is a chart of $\N$
by Lemma \ref{basic2}, Lemma \ref{basic3}, and (\ref{basic1}).
Take a chart
   $$
   \begin{array}{ccc}
   P&\longrightarrow&M\\
   \Bigdownarrow&&\Bigdownarrow\\
   \P&\longrightarrow&\M
   \end{array}
   $$
of $(X,\M)\rightarrow(Z,\P)$ extending $P\rightarrow\P$. Set $W=X\times_{Z}Y$.
There exists an induced homomorphism $M\rightarrow\O_W$. Define a log
structure on $W$ by this homomorphism. Then this
log scheme $(W,M\rightarrow\O_W)$ is the fiber product of (\ref{fpro1}) in
$\LS$. Since the associated log structure of $M\rightarrow\O_W$ is fine and
saturated, $(W,M)$ is, indeed, the fiber product of (\ref{fpro1}) in $\LSfs$.
\qed

\begin{dfn}\label{defthick}{\rm
The exact closed immersion $t:(T',\L')\rightarrow(T,\L)$
is said to be a {\it thickening of order $\leq n$}, if
$\I=\Ker(\O_T\rightarrow\O_{T'})$ is a nilpotent ideal such that
$\I^{n+1}=0$. }
\end{dfn}

\begin{lem}\label{thick}{\rm (cf.\ \cite {Ill1}).}
Let $(T,\L)$ and $(T',\L')$ be fine log schemes.
If $(t,\theta):(T',\L')\rightarrow(T,\L)$ is a thickening of order 1,
there exists a commutative diagram with exact rows:
   $$
   \begin{array}{ccccccccc}
   1&\rightarrow&1+\I&\hookrightarrow&t^{-1}\L&
   \stackrel{\theta}{\rightarrow}&\L'&\rightarrow&1\\
   &&\parallel&&\cap&&\cap\\
   1&\rightarrow&1+\I&\rightarrow&t^{-1}\gp{\L}&
   \underrel{\rightarrow}{\gp{\theta}}&\gp{\L'}&\rightarrow&1\rlap{,}
   \end{array}
   $$
where $\I=\Ker(\O_T\rightarrow\O_{T'})$, such that the right square of this
commutative diagram is cartesian.
\end{lem}

\noindent
The proof is straightforward.
Note that the multiplicative monoid $1+\I$ is identified with the additive
monoid $\I$ by $1+x\mapsto x$ since $\I^2=0$.

\begin{dfn}{\rm (cf.\ \cite [(3.3)]{Kat1})
Let $f:(X,\M)\rightarrow(Y,\N)$ be a morphism in $\LSf$. $f$ is said to be
{\it log smooth} if the following conditions are satisfied:
   \begin{enumerate}
   \item $f$ is locally of finite presentation,
   \item for any commutative diagram
      $$
      \begin{array}{ccc}
      (T',\L')&\stackrel{s'}{\longrightarrow}&(X,\M)\\
      \llap{$t$}\Bigdownarrow&&\Bigdownarrow\rlap{$f$}\\
      (T,\L)&\underrel{\longrightarrow}{s}&(Y,\N)
      \end{array}
      $$
   in $\LSf$,
   where $t$ is a thickening of order 1, there exists \'{e}tale locally
   a morphism $g:(T,\L)\rightarrow(X,\M)$ such that $s'=g\circ t$ and
   $s=f\circ g$.
   \end{enumerate}}
\end{dfn}

\noindent
The proofs of the following two propositions are straightforward and are
left to the reader.

\begin{pro}\label{usulisse}
Let $f:(X,\M)\rightarrow(Y,\N)$ be a morphism in $\LSf$.
If $f$ is strict,
then $f$ is log smooth if and only if $f$ is smooth in the usual sense.
\end{pro}

\begin{pro}\label{bextlisse}
For $(S,\L)\in\Obj(\LSf)$ and $(X,\M),(Y,\N)\in\Obj(\LSf_{(S,\L)})$,
let $f:(X,\M)\rightarrow(Y,\N)$ be a morphism in $\LSf_{(S,\L)}$.
Assume that $f$ is log smooth. If $(S',\L')$ is a log scheme over $(S,\L)$,
then the induced morphism
   $$
   (X,\M)\times_{(S,\L)}(S',\L')\rightarrow(Y,\N)\times_{(S,\L)}(S',\L')
   $$
is also log smooth.
\end{pro}
   \section{Toroidal characterization of log smoothness}
The following theorem is due to Kato \cite{Kat1},
and we prove it in \S \ref{prf1} for the reader's convenience.

\begin{thm}\label{lisse}{\rm (\cite [(3.5)]{Kat1})}
Let $f:(X,\M)\rightarrow(Y,\N)$ be a morphism in $\LSf$. and $Q\rightarrow\N$
a chart of $\N$, where $Q$ is a finitely generated integral
monoid. Then the following conditions are equivalent.
   \begin{description}
   \item[{\rm 1.}] $f$ is log smooth.
   \item[{\rm 2.}] There exists \'{e}tale locally a chart
   $(P\rightarrow\M,Q\rightarrow\N,Q\rightarrow P)$ of $f$ extending
   $Q\rightarrow\N$, where $P$ is a finitely generated integral
   monoid, such that
      \begin{description}
      \item[{\rm (a)}] $\Ker(\gp{Q}\rightarrow\gp{P})$ and the torsion part of
      $\Coker(\gp{Q}\rightarrow\gp{P})$ are finite groups of orders invertible
      on $X$,
      \item[{\rm (b)}] $X\rightarrow Y\times_{\Spec \Z[Q]}\Spec \Z[P]$ is
      smooth (in
      the usual sense).
      \end{description}
   \end{description}
Moreover, if $f:(X,\M)\rightarrow(Y,\N)$ is a log smooth
morphism in $\LSfs$ and
$Q\rightarrow\N$ is a chart of $\N$ such that $Q$ is finitely generated,
integral and saturated,
then there exists a chart
$(P\rightarrow\M,Q\rightarrow\N,Q\rightarrow P)$ of $f$ as above
with $P$ also saturated.
\end{thm}

\begin{rem}\label{lisserem}{\rm
The proof of Theorem \ref{lisse} in \S \ref{prf1} shows that we can require
in the condition 2. (a) that $\gp{Q}\rightarrow\gp{P}$ is injective without
changing the conclusion.
}
\end{rem}

We give some important examples of log smooth morphisms in the following.
Let $k$ be a field.

\begin{dfn}\label{logpt}{\rm
A log structure on $\Spec k$ is called a log structure of a {\it logarithmic
point} if it is equivalent (Definition \ref{logequiv})
to the associated log structure of
$\alpha:Q\rightarrow k$, where $Q$ is a monoid having no
invertible element other than 1 and $\alpha$ is a homomorphism defined by
   $$
   \alpha(x)=\left\{
      \begin{array}{ll}
      1&\mbox{if $x=1$,}\\
      0&\mbox{otherwise.}
      \end{array}
   \right.
   $$
Note that this log structure is equivalent to
$Q\oplus k^{\times}\rightarrow k$.
We denote the log scheme obtained in this way by $(\Spec k,Q)$.
The log scheme $(\Spec k,Q)$ is called a
{\it logarithmic point}.
Especially, if $Q=\Na$, the logarithmic point $(\Spec k,\Na)$ is said to be the
{\it standard log point}. }
\end{dfn}

\noindent
If $k$ is algebraically closed, any log structure on $\Spec k$ is equivalent
to a log structure of logarithmic point (cf.\ \cite {Ill1}).
Note that if we set
$Q=\{1\}$, then the log structure of the logarithmic point induced by $Q$ is
the trivial log structure (Example \ref{trilog}).

\begin{exa}\label{toroex}{\rm
Let $P$ be a submonoid of a group $M=\Z^d$ such that $\gp{P}=M$ and that
$P$ is
saturated. Let $Q$ be a submonoid of $P$, which is saturated but not
necessarily saturated in $P$.
We assume the following:
   \begin{enumerate}
   \item the monoid $Q$ has no invertible element other than 1,
   \item the order of the torsion part of $M/\gp{Q}$ is invertible
         in $k$.
   \end{enumerate}
Let $R=\Z[1/N]$ where $N$ is the order of the torsion part of
$M/\gp{Q}$.
The latter assumption implies, by Theorem \ref{lisse}, that
$(\Spec R[P],P)\rightarrow(\Spec R[Q],Q)$ (see Example \ref{canlog})
is log smooth. Define $\Spec k\rightarrow\Spec R[Q]$ by
$\alpha:Q\rightarrow k$ as in Definition \ref{logpt}. Let $X$ be
a scheme over $k$ which is smooth over $\Spec k\times_{\Spec R[Q]}\Spec R[P]$.
Then we have a diagram
   $$
   \begin{array}{ccc}
   X\\
   \Bigdownarrow\\
   \Spec k\times_{\Spec R[Q]}\Spec R[P]&\longrightarrow&\Spec R[P]\\
   \Bigdownarrow&&\Bigdownarrow\\
   \Spec k&\longrightarrow&\Spec R[Q]\rlap{.}
   \end{array}
   $$
Define a log structure $\M$ on $X$ by the pull--back of the canonical log
structure on $\Spec R[P]$.
Then we have a morphism
   $$
   f:(X,\M)\longrightarrow(\Spec k,Q)
   $$
of fine saturated log schemes. This morphsim $f$
is log smooth by Proposition \ref{usulisse} and
Proposition \ref{bextlisse}.
We denote this log scheme $(X,\M)$ simply by $(X,P)$. }
\end{exa}

\begin{exa}\label{toroidal}{\rm (Toric varieties.)
In this and the following examples, we use the notation appearing in
Example \ref{toroex}. Let $\sigma$ be a cone in $N_{\R}=\R^d$ and
$\sigma^{\vee}$ be its dual cone in $M_{\R}=\R^d$. Set
$P=M\cap\sigma^{\vee}$ and $Q=\{0\}\subset P$. Then,
$\Spec k\times_{\Spec \Z[Q]}\Spec \Z[P]$ is $k$--isomorphic to $\Spec k[P]$
which is nothing but an affine toric variety.
Let $X\rightarrow\Spec k[P]$ be a smooth morphism.
Then $(X,P)\rightarrow\Spec k$ is log smooth. }
\end{exa}

\begin{exa}\label{ssreduc}{\rm (Variety with normal crossings.)
Let $\sigma$ be the cone in $M_{\R}=\R^d$ generated by $e_1,\ldots,e_d$, where
$e_i=(0,\ldots,0,1,0,\ldots 0)\,\mbox{($1$ at the $i$--th entry)}$,
$1\leq i\leq d$. Let $\tau$ be the subcone generated by
$a_1e_1+\cdots+a_de_d$ with positive integers $a_j$ for $j=1,\ldots,d$.
We assume that $\mbox{\rm GCD}(a_1,\ldots,a_d)(=N)$
is invertible in $k$. Set $R=\Z[1/N]$.
Then, by setting $P=M\cap\sigma$ and $Q=M\cap\tau$, we see that
$\Spec k\times_{\Spec R[Q]}\Spec R[P]$ is $k$--isomorphic to
$\Spec k[z_1,\ldots,z_d]/(z_1^{a_1}\cdots z_d^{a_d})$ and $f$ is induced by
   $$
   \begin{array}{ccl}
   \Na^d&\longrightarrow&k[z_1,\ldots,z_d]/(z_1^{a_1}\cdots z_d^{a_d})\\
   \llap{$\varphi$}\Biguparrow&&\Biguparrow\\
   \Na&\longrightarrow&k,
   \end{array}
   $$
where the morphism in the first row is defined by $e_i\mapsto z_i,\,
(1\leq i\leq d)$, and $\varphi$ is defined by
$\varphi(1)=a_1e_1+\cdots+a_de_d$.
Let $X\rightarrow \Spec k[z_1,\ldots,z_d]/(z_1^{a_1}\cdots z_d^{a_d})$
be a smooth morphism.
Then, $(X,\Na^d)\rightarrow(\Spec k,\Na)$ is log smooth. }
\end{exa}

The following theorem is an application of Theorem \ref{lisse}.
We prove it in \S \ref{prf2}.

\begin{thm}\label{toroch}
Let $X$ be an algebraic scheme over a field $k$, and $\M\rightarrow\O_X$ a
fine saturated
log structure on $X$. Then, the log scheme $(X,\M)$ is log smooth over
$\Spec k$ with trivial log structure if and only if there exist an open
\'{e}tale covering $\U=\{U_i\}_{i\in I}$ of $X$
and a divisor $D$ of $X$
such that:
\begin{description}
\item[{\rm 1.}] there exists a smooth morphism
$$
h_i:U_i\longrightarrow V_i
$$
where $V_i$ is a affine toric variety over $k$ for each $i\in I$,
\item[{\rm 2.}] the divisor $U_i\cap D$ of $U_i$ is the pull--back of the
union of the closure of
codimension 1 torus orbits of $V_i$ by $h_i$ for each $i\in I$,
\item[{\rm 3.}] the log structure $\M\rightarrow\O_X$ is equivalent to the
log structure $\O_X\cap j_{\ast}\O^{\times}_{X-D}\hookrightarrow\O_X$ where
$j:X-D\hookrightarrow X$ is the inclusion.
\end{description}
\end{thm}

\begin{cor}\label{toroch1}
Let $X$ be a smooth algebraic variety over a field $k$, and $\M\rightarrow
\O_X$ a fine saturated log structure on $X$. Then, the log scheme
$(X,\M)$ is log smooth over $\Spec k$ with trivial log structure if and only
if there exists a reduced normal crossing divisor $D$ of $X$ such that
the log structure $\M\rightarrow\O_X$ is equivalent to the log structure
$\O_X\cap j_{\ast}\O^{\times}_{X-D}\hookrightarrow\O_X$ where
$j:X-D\hookrightarrow X$ is the inclusion.
\end{cor}

   \section{Log differentials and log derivations}
In this section, we are going to discuss the log differentials and
log derivations. These objects are closely related with log smoothness,
and play important roles in the sequel.
To begin with, we introduce a
useful notation which we often use in the sequel. Let $(X,\M)$ be a
log scheme. If we like to omit writing the log structure $\M$, we write
this log scheme by $\X$ to distinguish from the underlying scheme $X$.

\begin{dfn}\label{logder}
{\rm (cf.\ \cite {Kat1}, in different notation)
Let $\X=(X,\M)$ and $\Y=(Y,\N)$ be fine log schemes,
and $f=(f,\varphi):\X\rightarrow\Y$ a morphism, where
$\varphi:f^{-1}\N\rightarrow\M$ is a homomorphism of sheaves of monoids.
   \begin{enumerate}
   \item Let ${\cal E}$ be an $\O_X$--module.
         The sheaf of {\it log derivations} $\Ders_{\Y}(\X,{\cal E})$
         of $\X$ to ${\cal E}$ over $\Y$
         is the sheaf of germs of couples $(D,\Dlog )$, where
         $D\in\Ders_Y(X,{\cal E})$ and $\Dlog :\M\rightarrow{\cal E}$,
         such that the following conditions are satisfied:
         \begin{description}
      \item[{\rm (a)}] $\Dlog (ab)=\Dlog (a)+\Dlog (b),\mbox{ for}\,a,b\in\M$,
      \item[{\rm (b)}] $\alpha(a)\Dlog (a)=D(\alpha(a)),\mbox{ for}\,a\in\M$,
      \item[{\rm (c)}] $\Dlog (\varphi(c))=0,\mbox{ for}\,c\in f^{-1}\N$.
         \end{description}
   \item The sheaf of {\it log differentials} of $\X$ over $\Y$ is the
         $\O_X$--module defined by
         $$
         \Omega^1_{\X/\Y}=
         [\Omega^1_{X/Y}\oplus(\O_X\otimes_{\Z}\gp{\M})]/{\cal K},
         $$
         where ${\cal K}$ is the $\O_X$--submodule generated by
         $$
         (d\alpha(a),0)-(0,\alpha(a)\otimes a)\;\mbox{and}\;
         (0,1\otimes\varphi(b)),
         $$
         for all $a\in\M$, $b\in f^{-1}\N$.
   \end{enumerate}}
\end{dfn}

\noindent
These are coherent $\O_X$--modules if $Y$ is locally noetherian and $X$
locally of finite type over $Y$ (cf.\ \cite {Ill1}).
The proofs of the following three propositions are found in
\cite [\S 3]{Kat1}.

\begin{pro}
Let $\X$, $\Y$, $f$, and ${\cal E}$ be the same as in Definition \ref{logder}.
Then there is a natural isomorphism
   $$
   \Homs_{\O_X}(\Omega^1_{\X/\Y},{\cal E})
   \stackrel{\sim}{\longrightarrow}
   \Ders_{\Y}(\X,{\cal E}),
   $$
   by $u\mapsto(u\circ d,u\circ \dlog )$,
   where $d$ and $\dlog$ are defined by
   $$
   d:\O_X\rightarrow\Omega^1_{X/Y}\rightarrow\Omega^1_{\X/\Y}
   $$
   and
   $$
   \dlog :\M\rightarrow\O_X\otimes_{\Z}\gp{\M}\rightarrow\Omega^1_{\X/\Y}.
   $$
\end{pro}

\begin{pro}\label{genbun}
Let $\X\stackrel{f}{\rightarrow}\Y\stackrel{g}{\rightarrow}\Zu$ be morphisms
of fine log schemes.
\begin{enumerate}
   \item[{\rm 1.}] There exists an exact sequence
   $$
   f^{*}\Omega^1_{\Y/\Zu}\rightarrow\Omega^1_{\X/\Zu}
   \rightarrow\Omega^1_{\X/\Y}\rightarrow 0.
   $$
   \item[{\rm 2.}] If $f$ is log smooth, then
   \begin{equation}\label{diffseq}
   0\rightarrow f^{*}\Omega^1_{\Y/\Zu}\rightarrow\Omega^1_{\X/\Zu}
   \rightarrow\Omega^1_{\X/\Y}\rightarrow 0
   \end{equation}
   is exact.
   \item[{\rm 3.}] If $g\circ f$ is log smooth and {\rm (\ref{diffseq})}
   is exact and splits locally, then $f$ is log smooth.
\end{enumerate}
\end{pro}

\begin{pro}\label{bungen}
If $f:\X\rightarrow\Y$ is log smooth, then $\Omega^1_{\X/\Y}$ is a locally
free $\O_X$--module of finite type.
\end{pro}

\begin{exa}\label{fan1}{\rm (cf.\ \cite [Chap. 3, \S(3.1)]{Oda1})
Let $X_{\Sigma}$ be a toric variety over a field $k$
 determined by a fan $\Sigma$ on $N_{\R}$ with
$N=\Z^d$. Consider the log scheme $(X_{\Sigma},\Sigma)$
(Example \ref{torlog}) over $\Spec k$.
Then we have isomorphisms of $O_X$--modules
   $$
   \Ders_k(\X,\O_X)\cong\O_X\otimes_{\Z}N\; \; {\rm and}\; \;
   \Omega^1_{\X/k}\cong\O_X\otimes_{\Z}M,
   $$
where $M=\Hom_{\Z}(N,\Z)$. }
\end{exa}

\begin{exa}{\rm
For $X=\Spec k[z_1,\ldots,z_{n}]/(z_1\cdots z_l)$,
let $f:(X,\M)\rightarrow(\Spec k,\Na\rightarrow k)$ be the log smooth
morphism defined in Example \ref{ssreduc}. Then $\Ders_{\k}(\X,\O_X)$
is a free $\O_X$--module generated by
   $$
   z_1\frac{\partial}{\partial z_1},\ldots,z_l\frac{\partial}{\partial z_l},
   \frac{\partial}{\partial z_{l+1}},\ldots,\frac{\partial}{\partial z_n}
   $$
with a relation
   $$
   z_1\frac{\partial}{\partial z_1}+\cdots+z_l\frac{\partial}{\partial z_l}=0.
   $$
The sheaf $\Omega^1_{\X/\k}$ is a free $\O_X$--module generated by the
{\it logarithmic differentials}:
   $$
   \frac{dz_1}{z_1},\ldots,\frac{dz_l}{z_l},dz_{l+1},\ldots,dz_n
   $$
with a relation
   $$
   \frac{dz_1}{z_1}+\cdots+\frac{dz_l}{z_l}=0.
   $$
In the complex analytic case, the sheaf $\Omega^1_{\X/\k}$ is nothing but the
sheaf of {\it relative logarithmic differentials} introduced in, for
example, \cite [\S 3]{Fri1}, and \cite [\S 2]{K-N1}. }
\end{exa}
\section{The proof of Theorem 4.1}\label{prf1}
In this section, we give a proof of Theorem \ref{lisse} due to
Kato \cite {Kat1}.
Before proving the general case,
we prove the following proposition.

\begin{pro}\label{canlisse}
Let $A$ be a commutative ring and $h:Q\rightarrow P$ a homomorphism of
finitely generated integral monoids. The homomorphism $h$ induces the
morphism of log schemes
   $$
   f:\X=(\Spec A[P],P)\longrightarrow\Y=(\Spec A[Q],Q).
   $$
We set $K=\Ker(\gp{h}:\gp{Q}\rightarrow\gp{P})$ and
$C=\Coker(\gp{h}:\gp{Q}\rightarrow\gp{P})$, and denote the torsion part of $C$
by $\tor{C}$. If both $K$ and $\tor{C}$ are finite groups of order invertible
in $A$, then $f$ is log smooth.
\end{pro}

\pf
Suppose we have a commutative diagram
   $$
   \begin{array}{ccccc}
   (T',\L')&\stackrel{s'}{\longrightarrow}&\X&=&(\Spec A[P],P)\\
   \llap{$t$}\Bigdownarrow&&\Bigdownarrow\rlap{$f$}\\
   (T,\L)&\underrel{\longrightarrow}{s}&\Y&=&(\Spec A[Q],Q)
   \end{array}
   $$
in $\LSf$, where the morphism $t$ is a thickening of order 1.
Since we may work \'{e}tale locally, we may assume that $T$ is affine.
Set
   $$
   \I=\Ker(\O_T\rightarrow \O_{T'}).
   $$
Since the morphism $t$ is a thickening of order 1, by Lemma \ref{thick},
we have the following commutative diagram with exact rows:
   $$
   \begin{array}{ccccccccc}
   1&\longrightarrow&1+\I&\longhookrightarrow&\L&\stackrel{t^*}
   {\longrightarrow}&
   \L'&\longrightarrow&1\\
   &&\parallel&&\cap&&\cap\\
   1&\longrightarrow&1+\I&\longrightarrow&\gp{\L}&
   \underrel{\longrightarrow}{\gp{(t^*)}}&\gp{\L'}&
   \longrightarrow&1\rlap{.}
   \end{array}
   $$

\noindent
Note that the right square of the above commutative diagram is cartesian.

\vspace{3mm}
First, consider the following commutative diagram with exact rows:
   $$
   \begin{array}{ccccccccccc}
   1&\longrightarrow&K&\longrightarrow&\gp{Q}
   &\stackrel{\gp{h}}{\longrightarrow}&\gp{P}&\longrightarrow&C
   &\longrightarrow&1\\
   &&\Bigdownarrow\rlap{$u$}&&\Bigdownarrow\rlap{$v$}&&
   \Bigdownarrow\rlap{$w$}\\
   1&\longrightarrow&1+\I&\longrightarrow&\gp{\L}
   &\underrel{\longrightarrow}{\gp{(t^*)}}&\gp{\L'}&\longrightarrow
   &1\rlap{.}
   \end{array}
   $$
The multiplicative monoid $1+\I$ is isomorphic to the additive monoid $\I$
by $1+x\mapsto x$ since $\I^2=0$. If the order of $K$ is invertible in $A$,
then we have $u=1$, and hence there exists a morphism $a':R\rightarrow
\gp{\L}$ with $R=\mbox{\rm Image}\,(\gp{h}:\gp{Q}\rightarrow\gp{P})$ such
that $a'\circ\gp{h}=v$ and $\gp{(t^*)}\circ a'=w$.

\vspace{3mm}
Next, we consider the following commutative diagram with exact rows:
   $$
   \begin{array}{ccccccccccc}
   &&1&\longrightarrow&R&\stackrel{i}{\longrightarrow}&\gp{P}
   &\longrightarrow&C&\longrightarrow&1\\
   &&&&\llap{$a'$}\Bigdownarrow&&\Bigdownarrow\rlap{$w$}\\
   0&\longrightarrow&\I&\longrightarrow&\gp{\L}
   &\underrel{\longrightarrow}{\gp{(t^*)}}&\gp{\L'}&\longrightarrow
   &1\rlap{.}
   \end{array}
   $$
We shall show that there exists a homomorphism $a'':\gp{P}\rightarrow\gp{\L}$
such that $a''\circ t=a'$ and $\gp{(t^*)}\circ a''=w$. The obstruction of
existence of $a''$ lies in $\mbox{\rm Ext}^1(C,\I)$. In general, if a
positive integer $n$ is invertible in $A$ then we have
$\mbox{\rm Ext}^1(\Z/n\Z,\I)=0$. Combining this
with $\mbox{\rm Ext}^1(\Z,\I)=0$,
we have $\mbox{\rm Ext}^1(C,\I)=0$ since the order of the torsion part of $C$
is invertible in $A$. Hence a homomorphism $a''$ exists.
Since the diagram
   $$
   \begin{array}{ccc}
   \L&\stackrel{t^*}{\longrightarrow}&\L'\\
   \cap&&\cap\\
   \gp{\L}&\underrel{\longrightarrow}{\gp{(t^*)}}&\gp{\L'}\rlap{.}
   \end{array}
   $$
is cartesian, we found a
homomorphism
   $$
   a:P\longrightarrow\L
   $$
such that $t^*\circ a=(s')^*$ and $a\circ h=s^*$.
Using this $a$, we can construct a
morphism of log schemes
   $$
   g:(T,\L)\longrightarrow\X=(\Spec A[P],P)
   $$
such that $g\circ t=s'$ and $s\circ g=f$. \qed

\vspace{3mm}
Now, let us prove Theorem 4.1. First, we prove the implication
$2\Rightarrow 1$. Let $R=\Z[1/(N_1\cdot N_2)]$ where $N_1$ is the order of
$\Ker(\gp{Q}\rightarrow\gp{P})$ and $N_2$ is the order of the torsion part
of $\Coker(\gp{Q}\rightarrow\gp{P})$. By the assumption (a), we have
$$
Y\times_{\Spec \Z[Q]}\Spec \Z[P]\cong Y\times_{\Spec R[Q]}\Spec R[P].
$$
Since $X\rightarrow Y\times_{\Spec R[Q]}\Spec R[P]$ is smooth by (b),
$f$ is log smooth due to Proposition \ref{usulisse},
Proposition \ref{bextlisse} and Proposition \ref{canlisse}.

\vspace{3mm}
Next, let us prove the converse. Assume the morphism
$f$ is log smooth. Then, the sheaf $\Omega^1_{\X/\Y}$ is a locally free
$\O_X$--module of finite type (Proposition \ref{bungen}). Take any point
$x\in X$. We denote by $\bar{x}$ a separable closure of $x$.

\begin{ste}{\rm
Consider the morphism of $\O_X$--modules
   $$
   1\otimes\dlog:\O_X\otimes_{\Z}\gp{\M}\longrightarrow
   \Omega^1_{\X/\Y},
   $$
which is surjective by the definition of $\Omega^1_{\X/\Y}$. Then we can take
elements $t_1,\ldots,t_r\in\M_{\bar{x}}$ such that the system
$\{\dlog t_i\}_{1\leq i\leq r}$ is a $\O_{X,\bar{x}}$--base of
$\Omega^1_{\X/\Y,\bar{x}}$.
Consider the homomorphism $\psi:\Na^r\rightarrow\M_{\bar{x}}$ defined by
   $$
   \Na^r\ni(n_1,\ldots,n_r)\mapsto t_1^{n_1}\cdots t_r^{n_r}\in\M_{\bar{x}}.
   $$
Combining this $\psi$ with the homomorphism $Q\rightarrow f^{-1}(\N)_{\bar{x}}
\rightarrow\M_{\bar{x}}$, we have a homomorphism $\varphi:H=\Na^r\oplus
Q\rightarrow\M_{\bar{x}}$.
}
\end{ste}

\begin{ste}{\rm
Let $k(\bar{x})$ denote the residue field at $\bar{x}$. We have a homomorphism
   \begin{equation}\label{bunchan}
   k(\bar{x})\otimes_{\Z}\Z^r\longrightarrow k(\bar{x})\otimes_{\Z}
   \Coker(f^{-1}(\gp{\N}/\Oi_Y)_{\bar{x}}\rightarrow\gp{\M}_{\bar{x}}/
   \Oi_{X,\bar{x}})
   \end{equation}
by $k(\bar{x})\otimes_{\Z}\gp{\psi}:k(\bar{x})\otimes_{\Z}\Z^r\rightarrow
k(\bar{x})\otimes_{\Z}\gp{\M}_{\bar{x}}$ and canonical projectiones
$\gp{\M}_{\bar{x}}\rightarrow\gp{\M}_{\bar{x}}/\Oi_{X,\bar{x}}\rightarrow
\Coker(f^{-1}(\gp{\N}/\Oi_Y)_{\bar{x}}\rightarrow\gp{\M}_{\bar{x}}/
\Oi_{X,\bar{x}})$.
We claim that this morphism (\ref{bunchan}) is surjective.
In fact, this morphism coincides with the composite morphism
   $$
   k(\bar{x})\otimes_{\Z}\Z^r\rightarrow k(\bar{x})\otimes_{\O_{X,\bar{x}}}
   \Omega^1_{\X/\Y,\bar{x}}\rightarrow k(\bar{x})\otimes_{\Z}
   \Coker(f^{-1}(\gp{\N}/\Oi_Y)_{\bar{x}}\rightarrow\gp{\M}_{\bar{x}}/
   \Oi_{X,\bar{x}}),
   $$
where the first morphism is induced by $\dlog\circ\psi$ and the second one
by the canonical projection, and these morphisms are clearly surjective.
Hence the morphism (\ref{bunchan}) is surjective.
On the other hand, the homomorphism
   $$
   \gp{Q}\longrightarrow f^{-1}(\N/\Oi_Y)_{\bar{x}}
   $$
is surjective since $Q\rightarrow\N$ is a chart of $\N$. Hence, the
homomorphism
   $$
   k(\bar{x})\otimes_{\Z}\gp{Q}\longrightarrow k(\bar{x})\otimes_{\Z}
   f^{-1}(\N/\Oi_Y)_{\bar{x}}
   $$
is surjective, and then, the homomorphism
   $$
   1\otimes_{\Z}\gp{\varphi}:k(\bar{x})\otimes_{Z}\gp{H}\longrightarrow
   k(\bar{x})\otimes_{\Z}(\M_{\bar{x}}/\Oi_{X,\bar{x}})
   $$
is surjective.
This shows that the cokernel $C=\Coker(\gp{\varphi}:\gp{H}\rightarrow
\gp{\M}_{\bar{x}}/\Oi_{X,\bar{x}})$ is annihilated by an integer $N$
invertible in $\O_{X,\bar{x}}$.
}
\end{ste}

\begin{ste}{\rm
Take elements $a_1,\ldots,a_d\in\gp{\M}_{\bar{x}}$ which generates $C$. Then
we can write $a_i^n=u_i\varphi(b_i)$ for $u_i\in\Oi_{X,\bar{x}}$ and
$b_i\in\gp{H}$, for $i=1,\ldots,d$. Since $\Oi_{X,\bar{x}}$ is
$N$--divisible, we can write $u_i=v_i^N$ for $v_i\in\Oi_{X,\bar{x}}$, for
$i=1,\ldots,d$, and hence we may suppose $a_i^N=\varphi(b_i)$,
replacing $a_i$ by $a_i/v_i$, for $i=1,\ldots,d$. Let $G$ be the push--out
of the diagram
   $$
   \gp{H}\longleftarrow\Z^d\longrightarrow\Z^d,
   $$
where $\Z^d\rightarrow\gp{H}$ is defined by $e_i\mapsto b_i$,
and $\Z^d\rightarrow\Z^d$ is defined by $e_i\mapsto
Ne_i$ for $i=1,\ldots,d$. Then $\gp{\varphi}:
\gp{H}\rightarrow\gp{\M}_{\bar{x}}$ and $\Z^d\rightarrow\gp{\M}_{\bar{x}}$,
defined by $e_i\mapsto a_i$ for $i=1,\ldots,d$, induce the homomorphism
   $$
   \phi:G\longrightarrow\gp{\M}_{\bar{x}}
   $$
which maps $G$ surjectively onto $\gp{\M}_{\bar{x}}/\Oi_{X,\bar{x}}$.
Define $P=\phi^{-1}(\M_{\bar{x}})$, then $P$ defines a chart of $\M$ on some
neighborhood of $\bar{x}$ (\cite [Lemma 2.10]{Kat1}).
If $\M$ is saturated, $P$ is also saturated.
There exists an induced map $Q\rightarrow P$ which defines a chart of $f$ on
some neighborhood of $\bar{x}$.
Since $\gp{H}\rightarrow\gp{P}$ is injective, $\gp{Q}\rightarrow\gp{P}$ is
injective. The cokernel $\Coker(\gp{H}\rightarrow\gp{P})$ is annihilated
by $N$, hence $\tor{\Coker(\gp{Q}\rightarrow\gp{P})}$ is finite and
annihilated by $N$.
}
\end{ste}

\begin{ste}{\rm
Set $X'=Y\times_{\Spec \Z[Q]}\Spec \Z[P]$ and $g:X\rightarrow X'$. We need to
show that the morphism $g$ is smooth in the usual sense.
Since $\X$ has the log structure induced by $g$ from $\Xp=(X',P)$, it suffice
to show that $g$ is log smooth (Proposition \ref{usulisse}). Since
$k(\bar{x})\otimes_{\Z}(\gp{P}/\gp{Q})\cong k(\bar{x})\otimes_{\Z}\Z^d\cong
k(\bar{x})\otimes_{\O_{X,\bar{x}}}\Omega^1_{\X/\Y,\bar{x}}$ and
$\Omega^1_{\X/\Y}$ is locally free, we have $\Omega^1_{\X/\Y}\cong
\O_X\otimes_{\Z}(\gp{P}/\gp{Q})$ on some neighborhood of $\bar{x}$.
On the other hand, by direct calculations, one sees that
$\Omega^1_{\Xp/\Y}\cong\O_{X'}\otimes_{\Z[P]}\Omega^1_{\Z[P]/\Z[Q]}\cong
\O_{X'}\otimes_{\Z}(\gp{P}/\gp{Q})$. Hence we have $g^*\Omega^1_{\Xp/\Y}
\cong\Omega^1_{\X/\Y}$. This implies $g$ is log smooth due to Proposition
\ref{genbun} (in fact, $g$ is {\it log \'{e}tale} (cf.\ \cite{Kat1})).
}
\end{ste}

\noindent
This completes the proof of the theorem. \qed
   \section{The proof of Theorem 4.7}\label{prf2}
In this section, we give a proof of Theorem \ref{toroch}.
If $V=\Spec k[P]$ is an affine toric variety, then it is easy to
see that the log structure associated to $P\rightarrow k[P]$ is equivalent
to the log structure $\O_X\cap j_{\ast}\O^{\times}_{V-D}\hookrightarrow
\O_X$ where $D$ is the union of the closure of
codimention 1 torus orbits of $V$ and
$j:V-D\hookrightarrow V$ is the inclusion. Hence, the ``if'' part of
Theorem \ref{toroch} is easy to see. Let us prove the converse.
Let $(X,\M)$ be as in the assumption of Theorem \ref{toroch} and
$f:(X,\M)\rightarrow\Spec k$ the structure morphism.
The key--lemma is the following.

\begin{lem}\label{keylem1}
We can take \'{e}tale locally a chart $P\rightarrow\M$ of $\M$
such that
\begin{description}
\item[{\rm 1.}] the chart $(P\rightarrow\M,1\rightarrow k^{\times},
1\rightarrow P)$ of $f$ satisfies the conditions {\rm (a)} and {\rm (b)} in
Theorem \ref{lisse},
\item[{\rm 2.}] $P$ is a finitely generated integral saturated monoid, and
has no torsion element.
\end{description}
Here, by a torsion element, we mean an element $x\neq 1$ such that
$x^n=1$ for some positive integer $n$.
\end{lem}

First, we are going to show that the theorem follows from the above lemma.
Since the monoid $P$ has no torsion element, $P$ is the saturated submonoid of
a finitely generated free abelian group $\gp{P}$. Hence, $X$ is \'{e}tale
locally smooth over affine toric varieties, and the log structure $\M$ on $X$
is \'{e}tale locally equivalent to the pull--back of the log structure
induced by the union of the closure of
codimension 1 torus orbits. Since, these log structure
glue to the log structure $\M$ on $X$, the pull--back of the union
of the closure of codimension 1
torus orbits glue to a divisor on $X$. In fact,
this divisor is the compliment of the
largest open subset $U$ such that $\M|_U$ is trivial with the reduced scheme
structure.
Hence our assertion is proved.

\vspace{3mm}
Now, we are going to prove Lemma \ref{keylem1}. We may work \'{e}tale locally.
Take a chart $(P\rightarrow\M,1\rightarrow k^{\times},1\rightarrow P)$ of
$f$ as in Theorem \ref{lisse}. We may assume that $P$ is saturated.
Define
   $$
   \tor{P}=\{x\in P\: |\: x^n=1\,(\mbox{for some $n$})\}.
   $$
$\tor{P}$ is a subgroup in $P$.
Take a decomposition $\gp{P}=\fr{G}\oplus\tor{G}$ of the finitely
generated abelian group $\gp{P}$, where
$\fr{G}$ (resp.\ $\tor{G}$) is a free (resp.\ torsion) subgroup of
$\gp{P}$. Then we have the equalities $\tor{P}=P\cap\tor{G}=\tor{G}$ since
$P$ is saturated. Define a submonoid $\fr{P}$ by $\fr{P}=P\cap\fr{G}$.

\begin{cla}
$P=\fr{P}\oplus\tor{P}$.
\end{cla}

\pf
Take $x\in P$. Decompose $x=yz$ in $\gp{P}$ such that $y\in\fr{G}$ and
$z\in\tor{G}=\tor{P}$. Since $y^n=(xz^{-1})^n=x^n\in P$ for a large $n$,
we have $y\in P$. Hence $y\in\fr{P}$.
\qed

\vspace{3mm}\noindent
Define $\fr{\alpha}:\fr{P}\rightarrow\O_X$ by
$\fr{P}\hookrightarrow P\stackrel{\alpha}{\rightarrow}\O_X$.

\begin{cla}
The homomorphism
$\fr{\alpha}:\fr{P}\rightarrow\O_X$ defines a log structure
equivalent to $\M$.
\end{cla}

\pf
If $x\in\tor{P}$, then $\alpha(x)\in\Oi_X$ since
$\alpha(x)^n=1$ for a large $n$.
Hence $\alpha(\tor{P})\subset\Oi_X$. This implies that the associated log
structure of $\fr{P}$ is equivalent to that of $P$. \qed

\vspace{3mm}\noindent
Hence, the morphism $f$ is equivalent to the morphsim induced by the diagram
   $$
   \begin{array}{ccc}
   \fr{P}&\stackrel{\fr{\alpha}}{\longrightarrow}&\O_X\\
   \llap{$\fr{\varphi}$}\Biguparrow&&\Biguparrow\\
   1&\underrel{\longrightarrow}{\lambda}&k,
   \end{array}
   $$
Then we have to check the conditions (a) and (b) in
Theorem \ref{lisse}. The condition (a) is easy to verify.
Let us check the condition (b). We need to show that the morphism
   $$
   X\longrightarrow\Spec k[\fr{P}]
   $$
induced by $X\rightarrow\Spec\Z[P]\rightarrow\Spec\Z[\fr{P}]$ is smooth.

\begin{cla}
The morphism
   \begin{equation}\label{sepext}
   \Spec k[P]\longrightarrow
   \Spec k[\fr{P}]
   \end{equation}
induced by $\fr{P}\hookrightarrow P$ is \'{e}tale.
\end{cla}

\pf
Since $P=\fr{P}\oplus\tor{P}$, we have
$k[P]=k[\fr{P}]\otimes_{k}k[\tor{P}]$.
Since every element
in $\tor{P}$ is roots of 1, and the order of $\tor{P}$ is
invertible in $k$, the morphism
   $$
   k\longhookrightarrow k[\tor{P}]
   $$
is a finite separable extension of the field $k$.
This shows that the morphism (\ref{sepext}) is \'{e}tale. \qed

\vspace{3mm}\noindent
Now we have proved Lemma \ref{keylem1}, and hence, Theorem \ref{toroch}.
   \section{Formulation of log smooth deformation}
{}From now on, we fix the following notation. Let $k$ be a field and
$Q$ a finitely generated integral saturated monoid having no
invertible element other than 1.
Then we have a logarithmic point (Definition \ref{logpt}) $\k=(\Spec k, Q)$.
Let $f=(f,\varphi):\X=(X,\M)\rightarrow\k=(\Spec k,Q)$ be
a log smooth morphism in $\LSfs$.

\vspace{3mm}
Let $\Lambda$ be a complete noetherian local ring with the residue field $k$.
For example, $\Lambda=k$ or $\Lambda=\mbox{the Witt vector ring of $k$}$
when $k$ is perfect.
We denote, by $\Lambda[[Q]]$, the completion of the monoid ring
$\Lambda[Q]$ along
the maximal ideal $\mu+\Lambda[Q\setminus\{1\}]$ where $\mu$ denotes the
maximal ideal of $\Lambda$. The completion $\Lambda[[Q]]$ is a complete local
$\Lambda$--algebra and is noetherian since $Q$ is finitely generated.
If the monoid $Q$ is isomorphic to $\Na$ then the ring $\Lambda[[Q]]$ is
isomorphic to $\Lambda[[t]]$ as local
$\Lambda$--algebras.
Let $\art$ be the category of artinian local $\Lambda[[Q]]$--algebras with
the residue field $k$, and $\Art$ be the category of pro--objects
of $\art$ (cf.\ \cite {Sch1}).
For $A\in\Obj(\Art)$, we define a log structure on the scheme $\Spec A$ by
the associated log structure
   $$
   Q\oplus A^{\times}\longrightarrow A
   $$
of the homomorphism $Q\rightarrow\Lambda[[Q]]\rightarrow A$.
We denote, by $(\Spec A, Q)$, the log scheme obtained in this way.

\begin{dfn}{\rm
For $A\in\Obj(\art)$, a {\it log smooth lifting} of
$f:(X,\M)\rightarrow(\Spec k,Q)$ on $A$ is a morphism
$\widetilde{f}:(\widetilde{X},\widetilde{\M})\rightarrow(\Spec A,Q)$
in $\LSfs$ together with a cartesian diagram
   $$
   \begin{array}{ccc}
   (X,\M)&\longrightarrow&(\Xt,\Mt)\\
   \llap{$f$}\Bigdownarrow&&\Bigdownarrow\rlap{$\ft$}\\
   (\Spec k,Q)&\longrightarrow&(\Spec A,Q)
   \end{array}
   $$
in $\LSfs$.
Two liftings are said to be {\it isomorphic} if they are isomorphic in
$\LSfs_{(\Spec A,Q)}$. }
\end{dfn}

\noindent
Note that $(\Spec k,Q)\rightarrow(\Spec A,Q)$
is an exact closed immersion, and hence,
the above diagram is cartesian in $\LSfs$ if and only if so is in $\LS$
(Lemma \ref{fpro}).
In particular, the underlying morphisms of log smooth liftings are
(not necessarily flat) liftings in the usual sense.
Moreover, since exact closed immersions are stable under base changes,
$(X,\M)\rightarrow(\Xt,\Mt)$ is also an exact closed immersion.
If either $Q=\{1\}$ or $Q=\Na$, the underlying morphisms of any
log smooth liftings of $f$ are flat since these morphisms of log schemes are
{\it integral} (cf.\ \cite {Kat1}).
Hence, in this case, the underlying morphisms
of log smooth liftings of $f$ are flat liftings of $f$.

\vspace{3mm}
Take a local chart $(P\rightarrow \M, Q\rightarrow Q\oplus k^{\times},
Q\rightarrow P)$ of $f$
extending the given $Q\rightarrow k$ as in Theorem \ref{lisse} such that
$\gp{Q}\rightarrow\gp{P}$ is injective (Remark \ref{lisserem}). Then, $f$
factors throught $\Spec k\times_{\Spec \Z[Q]}\Spec \Z[P]$ by the smooth
morphism $X\rightarrow\Spec k\times_{\Spec \Z[Q]}\Spec \Z[P]$ and the
natural projection \'{e}tale locally. For $A\in\Obj(\art)$, a smooth lifting
   \begin{equation}\label{loclift}
   \Xt\longrightarrow\Spec A\times_{\Spec \Z[Q]}\Spec \Z[P]
   \end{equation}
of $X\rightarrow\Spec k\times_{\Spec \Z[Q]}\Spec \Z[P]$, with the naturally
induced log structure, gives a local log smooth lifting of $f$. Note that
this local lifting $(\Xt,\Mt)\rightarrow
(\Spec A,Q)$ is log smooth.
Conversely, suppose $\ft:(\Xt,\Mt)\rightarrow
(\Spec A,Q)$ is a local log smooth lifting of $f$ on $A$.

\begin{lem}\label{liftlem1}
The local chart $(P\rightarrow \M, Q\rightarrow Q\oplus k^{\times},
Q\rightarrow P)$ of $f$
lifts to the local chart
$(P\rightarrow\Mt, Q\rightarrow Q\oplus A^{\times}, Q\rightarrow P)$ of $\ft$.
\end{lem}

\pf
The proof is done by the induction with respect to the length of
$A$. Take $A'\in\Obj(\art)$ with surjective morphism
$A\rightarrow A'$ such that $I=\Ker(A\rightarrow A')\neq 0$ and $I^2=0$.
Let $\ft':(\Xt',\Mt')\rightarrow
(\Spec A',Q)$ be a pull--back of $\ft$. Then, $\ft'$ is a
log smooth lifting of $f$ to $A'$. By the induction, we have the lifted
local chart $(P\rightarrow\Mt', Q\rightarrow Q\oplus A'^{\times},
Q\rightarrow P)$
of $\ft'$. Since $(\Xt',\Mt')\rightarrow
(\Xt,\Mt)$ is a thickening of order 1, by Lemma
\ref{thick}, we have the
following commutative diagram with exact rows:
   $$
   \begin{array}{ccccccccc}
   0&\rightarrow&\I&\hookrightarrow&\Mt&
   \rightarrow&\Mt'&\rightarrow&1\\
   &&\parallel&&\cap&&\cap\\
   0&\rightarrow&\I&\rightarrow&\gp{\Mt}&
   \rightarrow&\gp{\Mt'}&\rightarrow&1\rlap{,}
   \end{array}
   $$
where $\I=\Ker(\O_{\Xt}\rightarrow\O_{\Xt'})$.
The right square of this diagram is cartesian. Consider the
following commutative diagram of abelian groups with exact rows and columns:
   $$
   \begin{array}{ccccccc}
   \Hom(\gp{P},\I)&\rightarrow&\Hom(\gp{P},\gp{\Mt})&\rightarrow&
   \Hom(\gp{P},\gp{\Mt'})&\rightarrow&\Extr^1(\gp{P},\I)\\
   \Bigdownarrow&&\Bigdownarrow&&\Bigdownarrow\\
   \Hom(\gp{Q},\I)&\rightarrow&\Hom(\gp{Q},\gp{\Mt})&\rightarrow&
   \Hom(\gp{Q},\gp{\Mt'})\\
   \Bigdownarrow\\
   \Extr^1(C,\I)\rlap{,}
   \end{array}
   $$
where $C=\Coker(\gp{Q}\hookrightarrow\gp{P})$. By our assumption, we have
$\Extr^1(\gp{P},\I)=0$ and $\Extr^1(C,\I)=0$, since the order of the
torsion part of each $\gp{P}$ and $C$ is invertible in $A$.
Then, by an easy diagram
chasing, we can show that the given morphism $\gp{P}\rightarrow
\gp{\Mt'}$ can be lifted to a morphism $\gp{P}\rightarrow
\gp{\Mt}$ such that
   $$
   \begin{array}{ccc}
   \gp{P}&\longrightarrow&\gp{\Mt}\\
   \Biguparrow&&\Biguparrow\\
   \gp{Q}&\longrightarrow&\gp{Q}\oplus A^{\times}
   \end{array}
   $$
is commutative. Then, we have the morphism $P\rightarrow\Mt$ and
we see that the diagram
   $$
   \begin{array}{ccc}
   \P&\longrightarrow&\Mt\\
   \Biguparrow&&\Biguparrow\\
   Q&\longrightarrow&Q\oplus A^{\times}
   \end{array}
   $$
is commutative. Since $\Mt/\Oi_{\Xt}
\stackrel{\sim}{\rightarrow}\Mt'/\Oi_{\Xt'}$, we can
easily show that the morphism $P\rightarrow\Mt$ defines a chart
(Lemma \ref {basic3}).
\qed

\vspace{3mm}
\noindent
Then, $\ft$ is factors throught $\Spec A\times_{\Spec \Z[Q]}\Spec \Z[P]$
by the induced morphism
$\Xt\rightarrow\Spec A\times_{\Spec \Z[Q]}
\Spec \Z[P]$ and the natural projection,
and we have the following commutative diagram
   $$
   \begin{array}{ccc}
   X&\longrightarrow&\Xt\\
   \Bigdownarrow&&\Bigdownarrow\\
   \Spec k\times_{\Spec \Z[Q]}\Spec \Z[P]&\longrightarrow&
   \Spec A\times_{\Spec \Z[Q]}\Spec \Z[P]\\
   \Bigdownarrow&&\Bigdownarrow\\
   \Spec k&\longrightarrow&\Spec A\rlap{,}
   \end{array}
   $$
such that the each square is cartesian. Hence,
$\Xt\rightarrow\Spec A\times_{\Spec \Z[Q]}
\Spec \Z[P]$ is smooth since it is a smooth lifting of the smooth morphism
$X\rightarrow\Spec k\times_{\Spec \Z[Q]}\Spec \Z[P]$, and is unique
up to isomorphisms by the classical theory. Therefore, we
have proved the following proposition.

\begin{pro}\label{liftloc}{\rm (cf. \cite [(3.14)]{Kat1})}
For $A\in\art$, a log smooth lifting of
$f:(X,\M)\rightarrow(\Spec k,Q)$ on $A$
exists \'{e}tale locally, and is unique up to isomorphisms. In particular,
log smooth liftings are log smooth.
\end{pro}

Let $\ft:(\Xt,\Mt)\rightarrow(\Spec A,Q)$ be a log smooth lifting of $f$ to
$A$, and $u:A'\rightarrow A$ a surjective homomorphism in $\art$ such that
$I^2=0$ where $I=\Ker (u)$. Suppose $\ft':(\Xt',\Mt')\rightarrow(\Spec A',Q)$
is a log smooth lifting of $f$ to $A'$ which is also a lifting of $\ft$.
Let $(P\rightarrow\Mt', Q\rightarrow Q\oplus A'^{\times}, Q\rightarrow P)$ be
a local chart of $\ft'$ which is a lifting of $(P\rightarrow\M, Q\rightarrow
Q\oplus k^{\times}, Q\rightarrow P)$. Define a local chart $(P\rightarrow\Mt,
Q\rightarrow Q\oplus A^{\times}, Q\rightarrow P)$ of $\ft$ by $P\rightarrow\M'
\rightarrow\M$ and $Q\rightarrow Q\oplus A'^{\times}\rightarrow Q\oplus
A^{\times}$. An automorphism $\Theta:(\Xt',\Mt')\stackrel{\sim}{\rightarrow}
(\Xt',\Mt')$ over $(\Spec A',Q)$ which is identity on $(\Xt,\Mt)$ induces
an automorphism $\theta:\gp{\Mt'}\stackrel{\sim}{\rightarrow}\gp{\Mt'}$.
Consider the diagram
   $$
   \begin{array}{ccccccccc}
   &&&&\gp{P}&=&\gp{P}\\
   &&&&\llap{$\alpha'$}\Bigdownarrow&&\Bigdownarrow\rlap{$\alpha$}\\
   1&\rightarrow&1+\I&\rightarrow&\gp{\Mt'}&\rightarrow&\gp{\Mt}&
   \rightarrow&1.
   \end{array}
   $$
For $a\in\gp{P}$, the element
$\alpha'(a)\cdot[\theta\circ\alpha'(a)]^{-1}$ is in
$1+\I$. Then, we have a morphism $\Delta:\gp{P}\rightarrow\I=I\cdot\O_{\Xt'}
\cong I\otimes_{A}\O_{\Xt}$ by
$\Delta(a)=\alpha'(a)\cdot[\theta\circ\alpha'(a)]^{-1}-1$.
The morphism $\Delta$
lifts to the morphism $\Delta:\gp{P}/\gp{Q}\rightarrow I\otimes_{A}\O_{\Xt}$
and defines a morphism of $\O_{\Xt}$--modules
   $$
   \O_{\Xt}\otimes_{\Z}(\gp{P}/\gp{Q})\rightarrow I\otimes_{A}\O_{\Xt}.
   $$
Since $\Omega^1_{\Xtd/\A}\cong\O_{\Xt}\otimes_{\Z}(\gp{P}/\gp{Q})$ \'{e}tale
locally, this defines a local section of
   $$
   \Homs_{\O_{\Xt}}(\Omega^1_{\Xtd/\A},I\otimes_{A}\O_{\Xt})\cong
   \Ders(\Xtd,\O_{\Xt})\otimes_{A}I.
   $$
Conversely, for a local section
$(D,\Dlog)\in\Ders(\Xtd,\O_{\Xt})\otimes_{A}I$,
$D$ induces an automorphism of $O_{\Xt'}$ and $\Dlog$ induces an
automorphism of $\Mt'$, and then, indues an automorphism of $(\Xt',\Mt')$.
By this, applying the arguement in SGA I \cite{Gro1} Expos\'{e} 3, we get
the following proposition.

\begin{pro}\label{liftauto}{\rm (cf. \cite [(3.14)]{Kat1})}
Let $\ft:(\Xt,\Mt)\rightarrow(\Spec A,Q)$ be a
log smooth lifting of $f$ to $A$, and $u:A'\rightarrow A$ a surjective
homomorphism in $\art$ such that $I^2=0$ where $I=\Ker (u)$
(i.e., $(\Spec A,Q)\rightarrow(\Spec A',Q)$ is a thickening of order
$\leq 1$).
   \begin{enumerate}
   \item The sheaf of germs of lifting automorphisms of $\Xt$ to
         $A'$ is
         $$
         \Ders_{\A}(\Xtd,\O_{\Xt})\otimes_A I.
         $$
   \item The set of isomorphism classes of log smooth liftings of
         $\Xt$ to $A'$ is isomorphic to
         $$
         \mbox{\rm H}^1(\Xt,\Ders_{\A}(\Xtd,\O_{\Xt}))\otimes_A I.
         $$
   \item The lifting obstructions of $\Xt$ to $A'$ are in
         $$
         \mbox{\rm H}^2(\Xt,\Ders_{\A}(\Xtd,\O_{\Xt}))\otimes_A I.
         $$
   \end{enumerate}
\end{pro}

Define the {\it log smooth deformation functor} $\D=\D_{\X/\k}$ by
   $$
   \D_{\X/\k}(A)=\{\mbox{isomorphism class of log smooth lifting of
   $f:\X\rightarrow\k$ on $A$}\}
   $$
for $A\in\Obj(\art$).
This is a covariant functor from $\art$ to $\Ens$, the category of sets,
such that $\D_{\X/\k}(k)$
consists of one point.
We shall prove the following theorem in the next section.

\begin{thm}\label{hull}
If the underlying scheme $X$ is proper over $k$,
then the log deformation functor $\D_{\X/\k}$ has a
representable hull {\rm (cf.\ \cite {Sch1})}.
\end{thm}
   \section{The proof of Theorem 8.5}
We are going to prove Theorem \ref{hull} by checking M. Schlessinger's
criterion (\cite [Theorem 2.11]{Sch1}) for $\D$. Let
$u_1:A_1\rightarrow A_0$ and $u_2:A_2\rightarrow A_0$
be morphisms in $\art$. Consider the map
   \begin{equation}\label{desant}
   \D(A_1\times_{A_0}A_2)\longrightarrow\D(A_1)\times_{\D(A_0)}\D(A_2).
   \end{equation}
Then we shall check the following conditions.
   \begin{description}
   \item[(H1)] The map (\ref{desant}) is a surjection whenever
                $u_2:A_2\rightarrow A_0$ is a surjection.
   \item[(H2)] The map (\ref{desant}) is a bijection when $A_0=k$ and
                $A_2=k[\epsilon]$, where $k[\epsilon]=k[E]/(E^2)$.
   \item[(H3)] $\mbox{dim}_k(t_{\D})<\infty$,
                where $t_{\D}=\D(k[\epsilon])$.
   \end{description}
By Proposition \ref{liftauto}, we have
   $$
   t_{\D}\cong\mbox{\rm H}^1(X,\Ders_{\k}(\X,\O_X)).
   $$
Our assumption implies that $t_{\D}$ is finite dimensional since
$\Ders_{\k}(\X,\O_X)$ is a coherent $\O_X$--module.
Hence we need to check (H1) and (H2).
Set $B=A_1\times_{A_0}A_2$. Let $v_i:B\rightarrow A_i$ be the natural map
for $i=1,2$. We denote the morphisms of schemes associated
to $u_i$ and $v_i$ also by $u_i:\Spec A_0\rightarrow\Spec A_i$ and
$v_i:\Spec A_i\rightarrow\Spec B$ for $i=1,2$, respectively.

\vspace{3mm}
\noindent{\sc Proof of (H1).}\hspace{2mm}
Suppose the homomorphism $u_2:A_2\rightarrow A_0$ is surjective.
Take an element
$(\eta_1,\eta_2)\in\D(A_1)\times_{\D(A_0)}\D(A_2)$ where $\eta_i$ is
an isomorphism class of a log smooth lifting
$f_i:(X_i,\M_i)\rightarrow(\Spec A_i,Q)$ for each $i=1,2$.
The equality
$\D(u_1)(\eta_1)=\D(u_2)(\eta_2)(=\eta_0)$ implies that there exists an
isomorphism
$(u_2)^*(X_2,\M_2)\stackrel{\sim}{\rightarrow}(u_1)^*(X_1,\M_1)$
over $(\Spec A_0,Q)$.
Here, $(u_i)^*(X_i,\M_i)$ is the pull--back of $(X_i,\M_i)$ by
$u_i:\Spec A_0\rightarrow\Spec A_i$ for $i=1,2$.
Set $(X_0,\M_0)=(u_1)^*(X_1,\M_1)$. We denote the induced morphism of
log schemes $(X_0, \M_0)\rightarrow(X_i, \M_i)$ by $\uu_i$ for
$i=1,2$. Then we have the following commutative diagram:
   $$
   \begin{array}{ccccc}
   (X_1,\M_1)&\stackrel{\uu_1}{\longleftarrow}&(X_0,\M_0)&
   \stackrel{\uu_2}{\longrightarrow}&(X_2,\M_2)\\
   \llap{$f_1$}\Bigdownarrow&&\llap{$f_0$}\Bigdownarrow&&
   \Bigdownarrow\rlap{$f_2$}\\
   (\Spec A_1,Q)&\stackrel{u_1}{\longleftarrow}&(\Spec A_0,Q)&
   \stackrel{u_2}{\longrightarrow}&(\Spec A_2,Q).
   \end{array}
   $$
We have to find an element $\xi\in\D(B)$, which represents a lifting of $f$
to $B$, such that $\D(v_i)(\xi)=\eta_i$ for $i=1,2$. Consider a scheme
$Z=(|X|,\O_{X_1}\times_{\O_{X_0}}\O_{X_2})$ over $\Spec B$. Define a log
structure on $Z$ by the natural homomorphism
   $$
   \N=\M_1\times_{\M_0}\M_2\longrightarrow
   \O_Z=\O_{X_1}\times_{\O_{X_0}}\O_{X_2}.
   $$
It is easy to verify that this homomorphism is a log structure. Since the
diagram
   $$
   \begin{array}{ccc}
   \N&\longrightarrow&\O_Z\\
   \Biguparrow&&\Biguparrow\\
   \llap{$Q\cong$\hspace{1mm}}Q\times_{Q}Q&\longrightarrow&B
   \end{array}
   $$
is commutative, we have a morphism $g:(Z,\N)\rightarrow(\Spec B,Q)$ of log
schemes. By the construction, we have the morphism
$v_i:(X_i,\M_i)\rightarrow (Z,\N)$ for $i=1,2$ such that the diagram
   $$
   \begin{array}{ccc}
   (X_1,\M_1)&\stackrel{\vv_1}{\longrightarrow}&(Z,\N)\\
   \llap{$\uu_1$}\Biguparrow&&\Biguparrow\rlap{$\vv_2$}\\
   (X_0,\M_0)&\underrel{\longrightarrow}{\uu_2}&(X_2,\M_2)
   \end{array}
   $$
is commutative. Since $u_2:A_2\rightarrow A_0$ is surjective, the underlying
morphism $X_1\rightarrow Z$ of $\vv_1$ is a closed immersion in the
classical sense. We have to show that the morphism $\vv_1$ is an
exact closed immersion. Take a local chart $(P\rightarrow\M, Q\rightarrow
Q\oplus k^{\times},
Q\rightarrow P)$ of $f$ as in Theorem \ref{lisse} such that $\gp{Q}
\rightarrow\gp{P}$ is injective. By Lemma \ref{liftlem1}, this local chart
lifts to the local chart of $f_i$ for each $i=0,1,2$. Since
$u_2:A_2\rightarrow A_0$ is surjective, we have the following isomorphism
   $$
   \N/\Oi_Z\stackrel{\sim}{\rightarrow}
   (\M_1/\Oi_{X_1})\times_{(\M_0/\Oi_{X_0})}(\M_2/\Oi_{X_2}).
   $$
By this, one sees that $P\cong P\times_{P}P\rightarrow\N$ is a local chart
of $\N$. This shows that $(Z,\N)$ is a fine saturated log scheme, and
$\vv_1$ is an exact closed immersion.
Then, the exact closed immersion $(X,\M)\rightarrow (X_1,\M_1)\rightarrow
(Z,\N)$ gives $g$ a structure of log smooth lifting of $f$ to $(\Spec B,Q)$.
Hence, $g$ represents an element $\xi\in\D(B)$.
It is easy to verify that $\D(v_i)(\xi)=\eta_i$ for $i=1,2$ since the morphism
$f_1$,$f_2$ and $g$ have the common local chart.
Thus (H1) is now proved.
\qed

\vspace{3mm}
\noindent{\sc Proof of (H2).}\hspace{2mm}
We continue to use the same notation as above. First, we prepare the following
lemma.

\begin{lem}\label{basic4}
Let $g':(Z',\N')\rightarrow (\Spec B,Q)$ be a log smooth lifting of $f$ with
a commutative diagram
   $$
   \begin{array}{ccc}
   (X_1,\M_1)&\longrightarrow&(Z',\N')\\
   \llap{$\uu_1$}\Biguparrow&&\Biguparrow\\
   (X_0,\M_0)&\underrel{\longrightarrow}{\uu_2}&(X_2,\M_2)
   \end{array}
   $$
of liftings such that $(v_i)^{*}(Z',\N')\stackrel{\sim}{\leftarrow}
(X_i,\M_i)$ over $(\Spec A_i,Q)$ for $i=1,2$. Then, the natural morphism
$(Z,\N)\rightarrow(Z',\N')$ is an isomorphism.
\end{lem}

\pf
We may work \'{e}tale locally. By Lemma \ref{liftlem1}, the local chart
$(P\rightarrow\M, Q\rightarrow Q\oplus k^{\times}, Q\rightarrow P)$
of $f$ lifts to a local
chart of $g'$. Take a local chart $(P\rightarrow\N, Q\rightarrow Q\oplus
B^{\times},
Q\rightarrow P)$ of $g$ by $P\rightarrow\N'\rightarrow\N$.
Then, the schemes $Z$ and $Z'$ are smooth liftings of
$X\rightarrow\Spec k\times_{\Spec \Z[Q]}\Spec \Z[P]$ to
$\Spec B\times_{\Spec \Z[Q]}\Spec \Z[P]$. Hence, we have only to show that
the natural morphism $Z\rightarrow Z'$ of underlying schemes is an
isomorphism. But this follows from the classical theory
\cite [Corollary 3.6]{Sch1} since each $X_i$ is a smooth lifting of
$X\rightarrow\Spec k\times_{\Spec \Z[Q]}\Spec \Z[P]$ to
$\Spec A_i\times_{\Spec \Z[Q]}\Spec \Z[P]$ for $i=0,1,2$.
\qed

\vspace{3mm}
\noindent
Let $g':(Z',\N')\rightarrow (\Spec B,Q)$ be a log smooth lifting of $f$
which represents a class $\xi'\in\D(B)$. Suppose that the class $\xi'$ is
mapped to $(\eta_1,\eta_2)$ by (\ref {desant}). Then,
   $$
   (X_0,\M_0)\stackrel{\sim}{\rightarrow}
   (v_1\circ u_1)_{*}(Z',\N')\cong (v_2\circ u_2)_{*}(Z',\N')
   \stackrel{\sim}{\leftarrow}(X_0,\M_0)
   $$
defines an automorphism $\theta$ of the lifting $(X_0,\M_0)$.
If this automorphism
$\theta$ lifts to an automorphism $\theta'$ of the lifting $(X_1,\M_1)$
such that $\theta'\circ \uu_1=\uu_1\circ\theta'$, then,
replacing $(X_1,\M_1)\rightarrow(Z',\N')$ by
$(X_1,\M_1)\stackrel{\theta'}{\rightarrow}(X_1,\M_1)\rightarrow(Z',\N')$,
we have a commutative diagram as in Lemma \ref{basic4}, and then, we have
$\xi=\xi'$. Now if $A=k$ (so that $(X_0,\M_0)=(X,\M)$, $\theta={\rm id}$),
$\theta'$ exists.
Thus (H2) is proved.
\qed

\vspace{3mm}\noindent
Thus we have proved Theorem \ref{hull}.
   \section{Example 1: Log smooth deformation over trivial base}
\label{exam1}
As we have seen in Theorem \ref{toroch}, any log scheme $(X,\M)$ which is
log smooth over $\Spec k$ with the trivial log structure is smooth over
an affine torus embedding \'{e}tale locally.

\begin{exa}{\rm (Usual smooth deformations.)
Let $X$ be a smooth algebraic variety over a field $k$.
Then $X$ with the trivial
log structure is log smooth over $\Spec k$
(this is the case $D=0$ in Corollary \ref{toroch1}), and our log smooth
deformation of $X$ is nothing but the usual smooth deformation of $X$. }
\end{exa}

\begin{exa}{\rm (Generalized relative deformations.)
Let $X$ be an algebraic  variety over a field $k$.
Assume that there exists a fine saturated log
structure $\M$ on $X$ such that $f:(X,\M)\rightarrow\Spec k$ is log
smooth. Then, by Theorem \ref{toroch}, $X$ is covered by \'{e}tale open
sets which are smooth over affine toric varieties, and the log structure
$\M\rightarrow\O_X$ is equivalent to the log structure defined by
   \begin{equation}\label{logdiv}
   \M=j_{\ast}\Oi_{X-D}\cap\O_X
   \end{equation}
for some divisor $D$ of $X$,
where $j$ is the inclusion $X-D\hookrightarrow X$.
In this situation, our log smooth deformation
of $f$ is the deformation of the pair $(X,D)$.
If $X$ is smooth over $k$, then $D$ is a reduced normal crossing
divisor (Corollary \ref{toroch1}).
Assume $X$ is smooth over $k$, then we have the exact sequence
   $$
   0\rightarrow\Ders_k(\X,\O_X)\rightarrow\Ders_k(X,\O_X)(=\Theta_X)
   \rightarrow{\cal N}\rightarrow 0,
   $$
where ${\cal N}$ is an $\O_X$--module locally written by
   $$
   ({\cal N}_{D_1|X}\otimes\O_{D_1})\oplus\cdots\oplus
   ({\cal N}_{D_d|X}\otimes\O_{D_d}),
   $$
where $D_1,\ldots,D_d$ are local components of $D$, and ${\cal N}_{D_i|X}$
is the normal bundle of $D_i$ in $X$ for $i=1,\ldots,d$.
Then, we have an exact sequence
   $$
   \mbox{\rm H}^0(D,{\cal N})\rightarrow t_{\D}\rightarrow
   \mbox{\rm H}^1(X,\Theta_X)\rightarrow
   \mbox{\rm H}^1(D,{\cal N}).
   $$
In this sequence, $\mbox{\rm H}^0(D,{\cal N})$ is viewed as the set of
isomorphism classes of locally trivial deformations of $D$ in $X$, and
$\mbox{\rm H}^1(D,{\cal N})$ is viewed as a set of obstructions of
deformations
of $D$ in $X$. Hence this sequence explains
the relation between the log smooth
deformation and the usual smooth deformation. Note that, if $D$ is a
smooth divisor on $X$, the log smooth deformation is nothing but the
{\it relative deformation} of the pair $(X,D)$ studied by Makio
\cite {Mak1}. }
\end{exa}

\begin{exa}{\rm (Toric varieties.)
Let $X_{\Sigma}$ be a complete toric variety
over a field $k$ defined by a fan $\Sigma$ in $N_{\R}$,
and consider the log scheme $\X=(X_{\Sigma},\Sigma)$
(Example \ref{torlog}) over $\Spec k$. We have seen in Example \ref{fan1} that
   $$
   \Ders_k(\X,\O_X)\cong\O_X\otimes_{\Z}N,
   $$
and hence is a globally free $\O_X$--module. Since
$\mbox{\rm H}^1(X,\Ders_k(\X,\O_X))=0$, any toric varieties are
infinitesimally rigid with respect to our log smooth deformation.
Note that toric varieties without log structures are not necessarily rigid
with respect to the usual smooth deformation. }
\end{exa}
   \section{Example 2: Smoothings of normal crossing varieties}
\label{exam2}
Let $k$ be a field. A {\it normal crossing variety} over $k$ is a seperated,
connected, and geometrically reduced scheme $X$ of finite type over $k$ which
is covered by an \'{e}tale open covering $\{X_\la\}_{\la\in\La}$ such that
each $X_{\la}$ is isomorphic to $\Spec k[z_1,\ldots,z_n]/(z_1\cdots z_{d_\la})$
over $k$ where $n-1={\rm dim}_k X$. Let $(\Spec k,\Na)$ be a standard log point
(Definition \ref{logpt}).

\begin{dfn}{\rm (cf. \cite[(2.6)]{Kaj1})\
A log structure $\alpha:\M\rightarrow\O_X$ on a normal crossing variety $X$
is called a log structure of {\it semistable type} over $(\Spec k,\Na)$ if
the following conditions are satisfied:
\begin{enumerate}
   \item there exists an \'{e}tale open covering $\{X_\la\}_{\la\in\La}$ of $X$
         such that, for each $\la\in\La$, $X_\la$ is isomorphic to
         $\Spec k[z_1,\ldots,z_n]/(z_1\cdots z_{d_\la})$ over $k$ and the log
         structure $\M_\la=\M|_{X_\la}$ has a chart $\beta_{\la}:\Na^{d_\la}
         \rightarrow\M_{\la}$ such that $\alpha\circ\beta_{\la}(e_i)=z_i$ for
         $i=1,\ldots, d_{\la}$, where $e_i=(0,\ldots,0,1,0,\ldots,0)$ ($1$ at
         the $i$--th entry),
   \item there exists a morphism $f:(X,\M)\rightarrow (\Spec k,\Na)$ which has
         a local chart $(\beta_{\la}:\Na^{d_\la}\rightarrow\M,
         \Na\rightarrow\Na\oplus
         k^{\times}, \varphi_{\la}:\Na\rightarrow\Na^{d_{\la}})$ where
         $\varphi_{\la}$ is the diagonal homomorphism, for each $\la\in\La$.
\end{enumerate}
}
\end{dfn}

\vspace{3mm}
\noindent
The morphism $f:(X,\M)\rightarrow (\Spec k,\Na)$ defined as above is called a
{\it logarithmic semistable reduction}. Note that a logarithmic semistable
reduction is log smooth (Example \ref{ssreduc}).
Let $f:(X,\M)\rightarrow(\Spec k,\Na)$ be a logarithmic semistable reduction.
Then, \'{e}tale locally, $X$ is isomorphic to
$k[z_1,\ldots,z_n]/(z_1\cdots z_d)$, and $f$ is induced by
the diagram
   $$
   \begin{array}{ccl}
   \Na^d&\longrightarrow&k[z_1,\ldots,z_n]/(z_1\cdots z_d)\\
   \llap{$\varphi$}\Biguparrow&&\Biguparrow\\
   \Na&\longrightarrow&k.
   \end{array}
   $$
For $A\in\Obj({\cal C}_{\Lambda[[\Na]]})$, the log structure on $A$ is
defined by $\gamma:\Na\rightarrow A$, $\gamma(1)=\pi\in m_A$, where $m_A$
is the maximal ideal of $A$. Then the log smooth lifting of $f$ on $A$
is locally equivalent to the morphism induced by the diagram
   $$
   \begin{array}{ccl}
   \Na^d&\longrightarrow&A[z_1,\ldots,z_n]/(z_1\cdots z_d-\pi)\\
   \llap{$\varphi$}\Biguparrow&&\Biguparrow\\
   \Na&\underrel{\longrightarrow}{\gamma}&A.
   \end{array}
   $$
Hence, our log smooth deformation carries out the smoothing of the normal
crossing variety $X$.

\vspace{3mm}
Now we assume that the singular locus $X_{\mbox{sing}}=D$ of
$X$ is connected. Then our deformation theory coincides with the
{\it logarithmic deformation theory of normal crossing varieties} introduced
by Kawamata and Namikawa \cite {K-N1} in the complex analytic
situation. The following theorem is essentially due to Kawamata and
Namikawa \cite {K-N1}, and we prove it in the next section.

\begin{thm}\label{dss}
Let $X$ be a normal crossing variety over $k$. Then, $X$ has a
log structure of semistable type over $(\Spec k,\Na)$
if and only if $X$ is $d$--semistable, i.e.,
${\cal E}xt^1(\Omega^1_X,\O_X)\cong\O_D$ {\rm (cf.\ \cite {Fri1})}.
\end{thm}

\section{The proof of Theorem 11.2}
Let $X$ be a normal crossing variety over a field $k$ of dimension $n$.
Then the scheme $X$ is covered by an \'{e}tale open covering
$\{X_{\la}\}_{\la\in\La}$ such that each $X_{\la}$ is isomorphic to the
divisor in the $(n+1)$--affine space over $k$, $V_\la=\Spec k[\zl_0,\ldots,
\zl_n]\cong \mbox{\bf A}^{(n+1)}_k$, defined by $\zl_0\cdots \zl_{d_\la}=0$
for some $0\leq d_\la \leq n$.
We call this covering $\cvr$ a {\it coordinate covering} of $X$.
If $\cvr$ is a coordinate covering, the {\it singular locus} of
$X$, denoted by $D$, is the closed subscheme of $X$ defined locally by
the equations $\zl_0\cdots\zl_{i-1}\cdot\zl_{i+1}\cdots\zl_{d_\la}=0$ for
$0\leq i \leq d_\la$ in $V_\la$.
We write the defining ideals of $X_\la$ and $D$
in $V_\la$ as
$$
  \begin{array}{ccl}
    \I_{X_\la} &=& (\zl_0\cdots \zl_{d_\lambda}),\vspace{2mm} \\
    \I_{D_\la} &=& (\zl_0\cdots\zl_{i-1}\cdot\zl_{i+1}\cdots
    \zl_{d_\la}\mid 0\leq i \leq d_\la).
  \end{array}
$$
These are ideals in $\O_{V_\la}$. In the sequel,
we fix these notation and conventions, sometimes omitting the
index $\la$.

\vspace{3mm}
Let us describe $\T^1_X=\Ext^1_{\O_X}(\Omega^1_X, \O_X)$
locally. Consider the following exact sequence,
$$
  0\rightarrow\I_X/\I^2_X\rightarrow\Omega^1_V\otimes\O_X
  \rightarrow\Omega^1_X\rightarrow 0.
$$
Here, we omitted the index $\la$. Its dual is
$$
  0\rightarrow\T^0_X\rightarrow\Theta_V\otimes\O_X
  \rightarrow\N_{X\mid V}.
$$
Then, we have an isomorphism $\T^1_X\cong
\Coker(\Theta_V\otimes\O_X\rightarrow\N_{X\mid V})$.

\begin{lem}\label{locT1}
  $\T^1_X$ is an invertible $\O_D$--module. More precisely,
  we have an isomorphism of $\O_D$--modules
  $$
    \T^1_X\cong(\I_X/\I^2_X)^\vee\otimes\O_D.
  $$
\end{lem}

\pf
  First, note that
  $\I_X/\I^2_X=\O_X\cdot d(z_0\cdots z_d)$, and
  $\Omega^1_V\otimes\O_X=\oplus^n_{i=1}\O_X\cdot dz_i$.
  Here, we omited the index $\la$. Let us denote
  the inclusion $\I_X/\I^2_X\hookrightarrow\Omega^1_V\otimes\O_X$ by
  $\iota$.
  For each $f\in\Theta_V\otimes\O_X=
  {\cal H}om_{\O_X}(\Omega^1_V\otimes\O_X, \O_X)$,
  set $f_i=f(dz_i)\in\O_X$. Then, $\iota^* f(d(z_0\cdots z_d))=
  \sum^d_{i=0}f_i z_0\cdots z_{i-1}\cdot z_{i+1}\cdots z_d\in\I_D\otimes\O_X$.
  Hence, $\Img\iota^*\subseteq{\cal H}om_{\O_X}(\I_X/\I^2_X, \I_D\otimes\O_X)$.
  Easy to see the converse. Then, we have $\Img\iota^*=
  {\cal H}om_{\O_X}(\I_X/\I^2_X, \I_D\otimes\O_X)$.
  Hence, we get
  $\T^1_X\cong\Coker\iota^*\cong{\cal H}om_{\O_X}(\I_X/\I^2_X, \O_D)
  \cong(\I_X/\I^2_X)^\vee\otimes\O_D$.
\qed

\vspace{3mm}
  Let $X$ be a reduced normal crossing variety over a field $k$
  and $\cvr$ a
  coordinate covering of $X$.
  A system $\{(\ztl_0,\ldots,\ztl_n)\}_{\la\in\La}$, where $\ztl_i\in
  {\rm H}^0(X_\la, \O_X)$ for $\la\in\La$ and $i=0,\ldots,n$,
  is said to be a {\it log system} on $X$ with respect to $\cvr$ if
  the following conditions are satisfied:
  \begin{enumerate}
    \item for $0\leq i\leq d_\la$, we have
          $
            \ztl_i=u^{\ssd{\la}}_i\cdot\zl_i
          $
          for some $u^{\ssd{\la}}_i\in{\rm H}^0(X_\la, \O^\times_X)$,
    \item for $d_\la <j\leq n$, $\ztl_j$ is invertible
          on $X_\la$.
  \end{enumerate}
Then we have the following lemma.

\begin{lem}\label{ltrns}
  On each $X_{\la\mu}=X_\la\cap X_\mu\neq\emptyset$,
  there exists a transition relation such as
  \begin{equation}\label{trans}
    \ztl_i=\ulm_i\ztm_{\slm(i)}\ (0\leq i\leq n),
  \end{equation}
  for some invertible function $\ulm_i$ on $X_{\la\mu}$ and a permutation
  $\slm$.
  {\rm (As is seen in the following proof, these $\ulm_i$'s
  and $\slm$'s are {\it not} unique. )}
\end{lem}

\pf
  Set $E^{\ssd{\la}}=\{i\mid(\ztl_i
  \mid_{X_{\la\mu}})\not\in{\rm H}^0(X_{\la\mu}, \O^\times_X)\}$,
  and set $E^{\ssd{\mu}}$ similarly.
  There is a one to one correspondence between elements in
  $E^{\ssd{\la}}$ and components of $X_{\la\mu}$ by
  $i\leftrightarrow\{\ztl_i=0\}$.
  This implies that there exists a bijection
  $\sigma: E^{\ssd{\la}}\stackrel{\sim}{\rightarrow}E^{\ssd{\mu}}$
  such that $\{\ztl_i=0\}=\{\ztm_{\sigma(i)}=0\}$.
  Hence, for $i\in E^{\ssd{\la}}$, we can write $\ztl_i
  =\ulm_i\cdot\ztm_{\sigma(i)}$ for some
  $\ulm_i\in{\rm H}^0(X_{\la\mu}, \O^\times_X)$.
  On the other hand, if $i$ is not in $E^{\ssd{\la}}$,
  we can write as (\ref{trans}) for these $i$'s,
  since $\ztl_i$ is invertible on $X_{\la\mu}$
  (note that the number of
  these $i$'s coincides with that of those $j$'s which are not in
  $E^{\ssd{\mu}}$).
\qed

\begin{pro}\label{cridss}{\rm (cf.\ \cite {K-N1})}
  The following conditions are equivalent.
  \begin{enumerate}
    \item $X$ is $d$--semistable.
    \item There exist a log system $\{(\ztl_0,\ldots,\ztl_n)\}$
    on $X$ and its transition system $\{(\ulm_i, \slm)\}$
    as in Lemma \ref{ltrns} such that the equality
    \begin{equation}\label{cdss}
      \ulm_0\cdots\ulm_n=1
    \end{equation}
    holds on each $X_\la\cap X_\mu\neq\emptyset$.
  \end{enumerate}
\end{pro}

\pf
  ($1\Rightarrow 2$) Assume that the normal crossing variety
  $X$ is $d$--semistable. Due to Lemma \ref{locT1}, the invertible
  $\O_X$--module
  $(\I_X/\I^2_X)\otimes\O_D\cong\I_X/\I_X\I_D$ is trivial, i.e.,
  $$
    \I_X/\I_X\I_D\quad\cong\quad\O_D.
  $$
  Let us denote the natural projection
  $\I_X/\I^2_X\rightarrow\I_X/\I_X\I_D$
  by $p$. The sheaf $\I_X/\I_X\I_D$ is a free $\O_D$-module and
  $p(z_0\cdots z_d)$ is a $\O_D$-free base of it. By the above isomorphism,
  the global section $1\in\O_D$ correspondes to $p(v\cdot z_0\cdots z_d)$,
  for some invertible function $v$ on $X$.
  Set
  $$
    \zeta_i=\left\{
    \begin{array}{ll}
      v\cdot z_0&(i=0), \\
      z_i&(1\leq i\leq d), \\
      1&(d<i\leq n).
    \end{array}
    \right.
  $$
  Then, the system $\{(\ztl_0,\ldots,\ztl_n)\}$ is a log system
  on $X$.
  Due to Lemma \ref{ltrns}, we can take a transition system
  $\{(\ulm_i, \slm)\}$ such that the transition relation
  (\ref{trans}) holds.
  On each $X_{\la\mu}\neq\emptyset$, we have $p_\la(\ztl_0\cdots\ztl_n)=
  p_\mu(\ztm_0\cdots\ztm_n)$.
  Then, we have $\ulm_0\cdots\ulm_n p_\mu(\ztm_0\cdots\ztm_n)=
  p_\mu(\ztm_0\cdots\ztm_n)$.
  Since, $p_\mu(v^{\ssd{\mu}}\zm_0\cdots\zm_{d_{\mu}})$
  is a $\O_D$--free base of $\I_X/\I_X\I_D$, we have
  $\ulm_0\cdots\ulm_n=1$ on $D$.
  Hence, we can write
  $$
    \ulm_0\cdots\ulm_n=1+\sum^n_{j=0}a_j
    \zm_0\cdots\zm_{j-1}\cdot\zm_{j+1}\cdots\zm_n.
  $$
  Replacing every $\ulm_j$ with
  $$
    \ulm_j-a_j\zm_0\cdots\zm_{j-1}\cdot\zm_{j+1}\cdots\zm_n
    (\ulm_0\cdots\ulm_{j-1}\cdot\ulm_{j+1}\cdots\ulm_n)^{-1},
  $$
  we get $\ulm_0\cdots\ulm_n=1$ on $X_{\la\mu}$ as desired.
  \vspace{3mm}

  ($2\Rightarrow 1$) The local section
  $\ztl_0\otimes\cdots\otimes\ztl_n$ is a local generator
  of $(\I_X/\I^2_X)\otimes\O_D$.
  By (\ref{cdss}), these local
  generators glue to a global section. Hence, the invertible
  $\O_D$--module $(\I_X/\I^2_X)\otimes\O_D$
  is trivial, and its dual $\T^1_X$ is also trivial.
\qed

\vspace{3mm}
Let $\cvr$ be a coordinate covering on $X$, and take a log system
$\{(\ztl_0,\ldots,\ztl_n)\}$ with respect to $\cvr$. On each $X_\la$,
consider a pre--log structure $\alpha_\la:\Na^{n+1}\longrightarrow\O_{X_\la}$,
defined by $\alpha_\la(e_i)=\ztl_i$ for $0\leq i\leq n$.
Then, the associate log structure of this pre--log structure is
equivalent to the log structure $\M_\la$ of semistable type
on the normal crossing variety $X_\la$.
We write this log structure, according to \S 1, by
\begin{equation}\label{loclog}
  \M_\la=\Na^{n+1}\oplus_{\alpha^{-1}_\la(\O^\times_{X_\la})}
  \O^\times_{X_\la}
  \stackrel{{\bar{\alpha}}_\la}{\longrightarrow}\O_{X_\la},
\end{equation}
and ${\bar{\alpha}}_\la(e_i, u)=u\cdot\ztl_i$,
for $0\leq i\leq n$.
Due to Lemma \ref{ltrns}, on each $X_{\la\mu}=X_\la\cap X_\mu\neq\emptyset$,
there exists a non canonical isomorphism
$$
  \begin{array}{ccc}
    \Na^{n+1}
    \oplus_{\alpha^{-1}_\la(\O^\times_{X_{\la\mu}})}
    \O^\times_{X_{\la\mu}} &
    \stackrel{\phi_{\la\mu}}{\longrightarrow} &
    \Na^{n+1}
    \oplus_{\alpha^{-1}_\mu(\O^\times_{X_{\la\mu}})}
    \O^\times_{X_{\la\mu}} \\
    \Bigdownarrow & & \Bigdownarrow \\
    \O_{X_{\la\mu}} & = & \O_{X_{\la\mu}}
  \end{array}
$$
defined by
\begin{equation}\label{isom}
  \phi_{\la\mu}(e_i, u)=(e_{\slm(i)}, u\cdot\ulm_i),
\end{equation}
for $0\leq i\leq n$.
\begin{lem}\label{patch}
  Any isomorphism $\phi_{\la\mu}$ which makes the above diagram commute
  is written in the form of (\ref{isom}).
\end{lem}

\pf
  Changing indicies suitably, we may assume that
  \begin{itemize}
    \item for $0\leq i\leq d$, $\alpha_\la(e_i)$ and
          $\alpha_\mu(e_i)$ are not invertible on $X_{\la\mu}$,
    \item for $d< i\leq n$, $\alpha_\la(e_i)$ and
          $\alpha_\mu(e_i)$ are invertible $X_{\la\mu}$.
  \end{itemize}
  Set $\phi_{\la\mu}(e_i, 1)=(\sum_j a_i^j e_j, \ulm_i)$,
  for $0\leq i\leq n$.
  \vspace{3mm}

  \noindent{\sc Case}\vspace{2mm} $0\leq i\leq d\:$:
  Since ${\bar{\alpha}}_\mu\circ\phi_{\la\mu}={\bar{\alpha}}_\la$,
  the matrix $(a_i^j)_{0\leq i,j\leq d}$ is a permutation matrix
  (Note that $A,B\in{\rm M}_d(\Na)$ and $AB=1$ implies that $A$ and $B$ are
  permutation matrices).
  Then we may assume $a_i^j=\delta_i^j\quad(0\leq i,j\leq d)$.
  Hence, we can write $\phi_{\la\mu}(e_i, 1)=(e_i, \ulm_i)\cdot(b_i, 1)$,
  where $\alpha_\mu(b_i)\in{\rm H}^0(X_{\la\mu}, \O^\times_X)$.
  Since we have the equality
  $$
    (e_i, \ulm_i\cdot\alpha_\mu(b_i))\cdot(b_i, \alpha_\mu(b_i)^{-1})
    =(e_i, \ulm_i)\cdot(b_i, 1),
  $$
  we get $\phi_{\la\mu}(e_i, 1)=(e_i, \uln_i\cdot\alpha_\mu(b_i))$
  in the quotient monoid
  $\mbox{\rm N}^{n+1}\oplus_{\alpha^{-1}_\mu(\O^\times_{X_{\la\mu}})}
  \O^\times_{X_{\la\mu}}$.
  \vspace{3mm}

  \noindent{\sc Case}\vspace{2mm} $d<i\leq n\:$: Since
  ${\bar{\alpha}}_\mu\circ\phi_{\la\mu}(e_i, 1)$ is invertible,
  we have $a_i^j=0$,
  for $d<i\leq n$ and $0\leq j\leq d$.
  This implies $\phi_{\la\mu}(e_i, 1)=(c_i, \ulm_i)$,
  where $\alpha_\mu(c_i)\in{\rm H}^0(X_{\la\mu}, \O^\times_X)$.
  Since
  $$
    (c_i, \ulm_i)\cdot(e_i, \alpha_\mu(e_i)^{-1})=
    (e_i, \alpha_\mu(c_i)\ulm_i\alpha_\mu(e_i)^{-1})\cdot
    (c_i, \alpha_\mu(c_i)^{-1}),
  $$
  we get $\phi_{\la\mu}(e_i, 1)=
  (e_i, \alpha_\mu(c_i)\cdot\uln_i\cdot\alpha_\mu(e_i)^{-1})$
  in the quotient monoid
  $\mbox{\rm N}^{n+1}\oplus_{\alpha^{-1}_\mu(\O^\times_{X_{\la\mu}})}
  \O^\times_{X_{\la\mu}}$.
  \vspace{3mm}

  Hence, for $0\leq i\leq n$, we have the equality $\phi_{\la\mu}(e_i, 1)=
  (e_i, \widetilde{\ulm_i})$,
  for some invertible function $\widetilde{\ulm_i}$ on $X_{\la\mu}$,
  in the quotient monoid
  $\mbox{\rm N}^{n+1}\oplus_{\alpha^{-1}_\mu(\O^\times_{X_{\la\mu}})}
  \O^\times_{X_{\la\mu}}$.
  Combining this with $\phi_{\la\mu}(0, u)=(0, u)$, we have the desired
  result.
\qed
\vspace{3mm}

\setcounter{ste}{0}
Now, let us prove Theorem 11.2.
  Assume that $X$ is $d$--semistable. Take a coordinate covering $\cvr$,
  a log system $\{(\ztl_0,\ldots,\ztl_n)\}$ and a transition system
  $\{(\ulm_i, \slm)\}$ such that (\ref{cdss}) holds.
  Then, each $X_\la$ has a log structure by (\ref{loclog}).
  Moreover, there exists a system of isomorphisms
  $\{\phi_{\la\mu}\}$ defined by (\ref{isom}). Now we are going to show
  that $\{\M_\la\}$ glues to a log structure $\M$ on $X$ of the desired type.

\begin{ste}{\rm
  Let us prove that
  $\phi_{\mu\nu}\circ\phi_{\la\mu}=\phi_{\la\nu}$ on each
  $X_{\la\mu\nu}=X_\la\cap X_\mu\cap X_\nu\neq\emptyset$, i.e.,
  $\{\phi_{\la\mu}\}$ satisfies the 1-cocycle condition.
  Set
  \begin{eqnarray*}
    \phi_{\mu\nu}\circ\phi_{\la\mu}(e_i, 1) & = &
    (e_{\tau\circ\sigma(i)}, v_{\sigma(i)}u_i) \\
    \phi_{\la\nu}(e_i, 1) & = & (e_{\rho(i)}, w_i),
  \end{eqnarray*}
  where we put $\sigma=\slm$, $\tau=\smn$, $\rho=\sln$,
  $u=\ulm$, $v=\umn$, and $w=\uln$.
  Changing indicies suitably, we may assume
  \begin{itemize}
    \item for $0\leq i\leq d$, $\ztl_i\mid_{X_{\la\mu\nu}}$
          is not invertible on $X_{\la\mu\nu}$,
    \item for $d<i\leq n$, $\ztl_i\mid_{X_{\la\mu\nu}}$
          is invertible on $X_{\la\mu\nu}$.
  \end{itemize}

  Since ${\bar{\alpha}}_\nu\circ\phi_{\mu\nu}\circ\phi_{\la\mu}=
  {\bar{\alpha}}_\la={\bar{\alpha}}_\nu\circ\phi_{\la\nu}$, we have
  \begin{equation}\label{str1}
    v_{\sigma(i)}u_i\ztn_{\tau\circ\sigma(i)}=\ztl_i=
    w_i\ztn_{\rho(i)}\quad(0\leq i\leq n).
  \end{equation}

  \vspace{3mm}
  \noindent
  {\sc Case}\hspace{2mm} $d<i\leq n\:$: Since the equality (\ref{str1})
  holds, both
  $\ztn_{\tau\circ\sigma(i)}$ and $\ztn_{\rho(i)}$ are invertible
  on $X_{\la\mu\nu}$ and we have
  $$
    (e_{\tau\circ\sigma(i)}, v_{\sigma(i)}u_i)\cdot
    (e_{\rho(i)}, (\ztn_{\rho(i)})^{-1})=
    (e_{\rho(i)}, w_i)\cdot
    (e_{\tau\circ\sigma(i)}, (\ztn_{\tau\circ\sigma(i)})^{-1}).
  $$
  Hence, we have $(e_{\tau\circ\sigma(i)}, v_{\sigma(i)}u_i)=
  (e_{\rho(i)}, w_i)$ in the quotient monoid
  $\Na^{n+1}
  \oplus_{\alpha^{-1}_\nu(\O^\times_{X_{\la\mu\nu}})}
  \O^\times_{X_{\la\mu\nu}}$.

  \vspace{3mm}
  \noindent
  {\sc Case}\hspace{2mm} $0\leq i\leq d\:$: Since the components
  $\{\ztn_{\tau\circ\sigma(i)}=0\}$ and $\{\ztn_{\rho(i)}=0\}$
  of $X_{\la\mu\nu}$ coincides due to (\ref{str1}),
  we have $\tau\circ\sigma(i)=\rho(i)$.
  Hence, we have $\ztn_{\rho(i)}=w_i^{-1}v_{\sigma(i)}u_i\ztn_{\rho(i)}$.
  This imples
  \begin{equation}\label{str2}
    w_i^{-1}v_{\sigma(i)}u_i=1+
    a_i\zn_{\rho(0)}\cdots\zn_{\rho(i-1)}\cdot\zn_{\rho(i+1)}
    \cdots\zn_{\rho(d)},
  \end{equation}
  for some $a_i\in\O_{X_{\la\mu\nu}}$.
  Since
  \begin{eqnarray*}
    \rho(\{d+1,\ldots,n\})&=&\tau\circ\sigma(\{d+1,\ldots,n\})\\
                          &=&\{j\mid(\ztn_j\mid_{X_{\la\mu\nu}})\in
                             {\rm H}^0(X_{\la\mu\nu}, \O^\times_X)\},
  \end{eqnarray*}
  and $\ztn_{\rho(j)}=w_j^{-1}v_{\sigma(j)}u_j
  \ztn_{\tau\circ\sigma(j)}$ for $d<j\leq n$, we have
  $$
    \prod_{d<j\leq n}w_j^{-1}v_{\sigma(j)}u_j=1.
  $$
  On the other hand, by our assumptions
  $u_1\cdots u_n=1$, etc., we have
  $$
    \prod_{0\leq i\leq n}w_j^{-1}v_{\sigma(j)}u_j=1.
  $$
  Hence, we get the equality
  $$
    \prod_{0\leq i\leq d}w_j^{-1}v_{\sigma(j)}u_j=1.
  $$
  By this and (\ref{str2}), we have
  $$
    \sum_{0\leq i\leq d}a_i\zn_{\rho(0)}\cdots
    \zn_{\rho(i-1)}\cdot\zn_{\rho(i+1)}\cdots
    \zn_{\rho(d)}=0,
  $$
  and consequently we get $a_i=0$, that is, $w_i^{-1}v_{\sigma(i)}u_i=1$.
  Hence, also in this case, we have
  $(e_{\tau\circ\sigma(i)}, v_{\sigma(i)}u_i)=(e_{\rho(i)}, w_i)$.
  Therefore, we have proved that $\phi_{\mu\nu}\circ\phi_{\la\mu}=
  \phi_{\la\nu}$ holds and $\{\M_\la\}$ glues to a log structure $\M$ on $X$.
}
\end{ste}

\begin{ste}{\rm The system of morphisms
  $\{f_\la\}$,
  where $f_\la:(X_\la, \M_\la)\rightarrow
  (\Spec k, \Na)$ is defined as in Definition \ref{canlog},
  glues to the morphism $f$ if and only if
  the following diagram commutes:
  $$
    \begin{array}{ccc}
      \Na^{n+1}\oplus_{\alpha^{-1}_\la(\O^\times_{X_{\la\mu}})}
      \O^\times_{X_{\la\mu}} & \stackrel{\phi_{\la\mu}}{\longrightarrow} &
      \Na^{n+1}\oplus_{\alpha^{-1}_\mu(\O^\times_{X_{\la\mu}})}
      \O^\times_{X_{\la\mu}} \\
      \Biguparrow & & \Biguparrow \\
      \Na\oplus k^\times
      & = & \Na\oplus k^\times,
    \end{array}
  $$
  where $\Na\oplus k^\times\rightarrow k$ is the associated log
  structure of $\Na\rightarrow k$, and the homomorphism
  $\Na\oplus k\rightarrow
  \Na^{n+1}\oplus_{\alpha^{-1}_\la
  (\O^\times_{X_{\la\mu}})}\O^\times_{X_{\la\mu}}$
  is defined by $(1, 1)\mapsto (e_0+\cdots+e_n, 1)$.
  The commutativity of this diagram is equivalent to the equality
  \begin{eqnarray*}
    (e_0+\cdots+e_n, 1) & = & \phi_{\la\mu}(e_0+\cdots+e_n, 1) \\
                        & = & (e_0+\cdots+e_n, \ulm_0\cdots\ulm_n).
  \end{eqnarray*}
  It is easy to see that this is equivalent to $\ulm_0\cdots\ulm_n=1$.
  Hence, $\{f_\la\}$ glues to the morphism
  $f$ due to Proposition \ref{cridss}.
  Thus, we have proved that $\M$ is the log structure on $X$ of the desired
  type.
}
\end{ste}

\begin{ste}{\rm Let us prove the converse.
  Assume that $X$ has a log structure of semistavle type.
  Then, by Lemma \ref{patch}, Step 2 above and
  Proposition \ref{cridss},
  it is easy to show that $X$ is $d$--semistable.
}
\end{ste}
Thus, the proof of Theorem 11.2 is completed.
%
%


\begin{thebibliography}{99}
\bibitem{Fri1}{\sc R. Friedman}, {\sl Global smoothings of varieties with
              normal crossings}, Annals of Math. {\bf 118} (1983), 75--114.
\bibitem{Gro1}{\sc A. Grothendieck}, {\sl Rev\^{e}tement \'{e}tale et groupe
              fondamental}, Springer Lecture Note {\bf 224} (1971).
\bibitem{Ill1}{\sc L. Illusie}, {\sl Introduction \`{a} la g\'{e}om\'{e}trie
              logarithmique}, Seminary note.
\bibitem{Kaj1}{\sc T. Kajiwara}, {\sl Logarithmic compactifications of the
              generalized Jacobian variety}, J. Fac. Sci. Univ. Tokyo,
              Sect. IA, Math. {\bf 40} (1993), 473--502.
\bibitem{Kat1}{\sc K. Kato}, {\sl Logarithmic structures of
              Fontaine--Illusie}, in Algebraic Analysis, Geometry and
              Number Theory, J.--I. Igusa ed., 1988, Johns Hopkins Univ.,
              191--224.
\bibitem{K-N1}{\sc Y. Kawamata and Y. Namikawa}, {\sl Logarithmic Deformations
              of Normal Crossing Varieties and Smoothings of Degenerate
              Calabi--Yau Varieties}, preprint, 1992.
\bibitem{K-S1}{\sc K. Kodaira and D.C. Spencer}, {\sl On deformation of
              complex analytic structures, I--II}, Annals of Math. {\bf 67}
              (1958), 328--466.
\bibitem{L-S1}{\sc S. Lichtenbaum and M. Schlessinger}, {\sl The Cotangent
              Complex of a Morphism}, Trans. A.M.S. {\bf 128} (1967), 41--70.
\bibitem{Mak1}{\sc K. Makio}, {\sl On the Relative Pseudo--Rigidity},
              Proc. Japan Acad., {\bf 49} (1973), 6--9.
\bibitem{Oda1}{\sc T. Oda}, {\sl Convex Bodies and Algebraic Geometry: An
              Introduction to the Theory of Toric Varieties}, Ergebnisse der
              Math. (3) {\bf 15}, Springer--Verlag, 1988.
\bibitem{Sch1}{\sc M. Schlessinger}, {\sl Functors of Artin rings},
              Trans. A.M.S. {\bf 130} (1968), 208--222.
\bibitem{Ste1}{\sc J. H. M. Steenbrink}, {\sl Logarithmic Embeddings of
              Varieties with Normal Crossings and Mixed Hodge Structures},
              Preprint, 1993.
\end{thebibliography}
\end{document}